\documentclass[final]{siamltex}

\usepackage{amsfonts}
\usepackage{psfrag,graphicx}

\newtheorem{remark}[theorem]{Remark}

\newcommand\sqr[2]{{\vcenter{\vbox{\hrule height.#2pt
   \hbox{\vrule width.#2pt height#1pt \kern#1pt
      \vrule width.#2pt}
   \hrule height.#2pt}}}}

\newcommand{\ds}{\displaystyle}

\renewcommand{\theequation}{\thesection.\arabic{equation}}

\renewcommand{\epsilon}{\varepsilon}
\newcommand{\eps}{\epsilon}
\renewcommand{\t}{\tau}

\newcommand{\BV}{\mathbf{BV}}

\newcommand{\C}[1]{{{\mathbf{C}^\mathbf{#1}}}}

\renewcommand{\L}[1]{{{\mathbf{L}^\mathbf{#1}}}}
\newcommand{\Lloc}[1]{\mathbf{L}^{\mathbf{1}}_{\mathrm{loc}}}

\newcommand{\reali}{{\mathbb{R}}}

\newcommand{\naturali}{\mathbb{N}}

\newcommand{\tv}{\mathrm{TV}}

\renewcommand{\sf}[2]{{B\!\left(#1,#2\right)}}
\newcommand{\piu}[1]{{{[\![\,{{#1}}\,]\!]_{\strut +}^{\vphantom{\strut +}}}}}
\newcommand{\meno}[1]{{{[\![\,{{#1}}\,]\!]_{\strut -}^{\vphantom{\strut +}}}}}

\makeatletter

\newlength{\captionwidth}
\long\def\@makecaption#1#2{%
   \vskip 10\p@
   \setbox\@tempboxa\hbox{#1: #2}%
   \ifdim \wd\@tempboxa > \captionwidth 
       \hbox to\hsize{\hfil
       \parbox[t]{\captionwidth}{
       #1: #2\par}
       \hfil}
     \else
       \hbox to\hsize{\hfil\box\@tempboxa\hfil}%
   \fi}

\makeatother

\newcommand{\sz}{\scriptsize}

\newcommand{\wtv}{\mathrm{WTV}}

\renewcommand{\piu}[1] {{{[{{#1}}]_{\strut +} }}}
\renewcommand{\meno}[1]{{{[{{#1}}]_{\strut -}}}}

\def\knp{K_{np}}

\title{On a model of multiphase flow}

\author{Debora Amadori\thanks{Dipartimento di Matematica Pura e Applicata,
Universit\`a degli Studi dell'Aquila, Via Vetoio, 67010 - Coppito
(AQ), Italy ({\tt amadori@univaq.it})}
\and Andrea Corli\thanks{Dipartimento di Matematica,
Universit\`a di Ferrara, Via Machiavelli 35, 44100 Ferrara,
Italy ({\tt andrea.corli@unife.it}). Please send any
correspondence to this author.}
}

\begin{document}

\maketitle

\begin{abstract}
  We consider a hyperbolic system of three conservation laws in one
  space variable. The system is a model for fluid flow allowing phase
  transitions; in this case the state variables are the specific
  volume, the velocity and the mass density fraction of the vapor in
  the fluid.  For a class of initial data having large total variation
  we prove the global existence of solutions to the Cauchy problem.
\end{abstract}






\section{Introduction}
We consider a model for the one-dimensional flow of an inviscid fluid
capable of undergoing phase transitions.  Both liquid and vapor phases
are possible, as well as mixtures of them. In Lagrangian coordinates
the model is

\begin{equation}\label{eq:system}
\left\{
\begin{array}{ll}
v_t - u_x &= 0\\
u_t + p(v,\lambda)_x &= 0
\\
\lambda_t &= 0\,.
\end{array}
\right.
\end{equation}
Here $t>0$ and $x\in\reali$; moreover $v>0$ is the specific
volume, $u$ the velocity, $\lambda$ the mass density fraction of
vapor in the fluid. Then $\lambda\in[0,1]$, with $\lambda=0$
characterizing the liquid and $\lambda=1$ the vapor phase; the
intermediate values of $\lambda$ model the mixtures of the two
pure phases. The pressure is denoted by $p=p(v,\lambda)$; under
natural assumptions the system is strictly hyperbolic.

This model is a simplified version of a model proposed by Fan
\cite{Fan}, where also viscous and relaxation terms were taken into
account. The model is isothermal, see (\ref{eq:pressure}) below; in
presence of phase transitions this physical assumption is meaningful
for retrograde fluids. A study of the Riemann problem for a $2 \times
2$ relaxation approximation of (\ref{eq:system}) has been done in
\cite{CorliFan}. We focus here on the global existence of solutions to
the Cauchy problem for (\ref{eq:system}), namely for initial data
$$
(v,u,\lambda)(0,x)=\bigl(v_o(x),u_o(x),\lambda_o(x)\bigr)
$$
having finite total variation. This problem is motivated by the
study of more complete models, where (\ref{eq:system}) is supplemented
by source terms.

\smallskip

The problem of the global existence of solutions to strictly
hyperbolic system of conservation laws has been studied since
long, see \cite{Bressanbook, Dafermosbook, Serre, Smoller} for
general information. If the initial data have \emph{small} total
variation then Glimm theorem \cite{Glimm} applies; we refer again
to \cite{Bressanbook} for the analogous results obtained by a
wave-front tracking algorithm as well as for uniqueness and
continuous dependence of the solutions on the initial data.

Some special systems allow however initial data with \emph{large}
total variation. For the system of isothermal gasdynamics Nishida
\cite{Nishida68} proved that it is sufficient that the variation
$\tv(v_o,u_o)$ of the initial data is finite in order to have
globally defined solutions. This result was extended by Nishida
and Smoller \cite{NishidaSmoller} to any pressure law
$p=k/v^\gamma$, $\gamma>1$, provided that $(\gamma-1)\tv(v_o,u_o)$
is small; related results are in \cite{DiPerna}. For the full
nonisentropic system of $3\times3$ gasdynamics,
$p=k\exp(\frac{\gamma-1}{R}s)/v^\gamma$, for $s$ the entropy, Liu
\cite{Liu, LiuIB} proved the global existence of solutions if
$(\gamma-1) \tv(v_o,p_o)$ is small and $\tv(s_o)$ bounded. Temple
\cite{Temple} and Peng \cite{Peng92} obtained similar results. All
these papers use the Glimm scheme. Analogous results making use of
a wave-front tracking scheme have been given recently by Asakura
\cite{Asakura2}, \cite{Asakura3}; we point out that the use of
wave-front tracking schemes in case of data with large variation
is far from being trivial, and a deep analysis of the wave
interactions is required. Very general results can be proved for
systems with coinciding shock and rarefaction curves,
\cite{Bianchini}; however system (\ref{eq:system}) is not of this
type.

\smallskip

In comparison with the above systems of gasdynamics, in
(\ref{eq:system}) we keep a $\gamma$-law for the pressure with
$\gamma=1$, but add a dependence of $p$ on $\lambda$: we take then
$p=a(\lambda)/v$ for a suitable function $a$. System (\ref{eq:system})
has close connections to a system introduced by Benzoni-Gavage
\cite{Benzoni-Gavage} and studied by Peng \cite{Peng}; it seems
however that the proof in \cite{Peng} is not complete. A comparison of
these models is done in Subsection \ref{sec:comparison}. We mention
that also the method of compensated compactness has been applied to
(\ref{eq:system}), see \cite{Gosse}, \cite{BBL} and \cite[\S 12.3, \S
16]{Lu}, but for different pressure laws.

\smallskip

In this paper we prove by a wave-front tracking scheme the global
existence of solutions to (\ref{eq:system}) for a wide class of
initial data with large total variation. We introduce first a
\emph{weighted total variation} ($\wtv$) 
of $a(\lambda_o)$; this quantity arises in a natural way in the
problem and has also an analytical meaning, being the logarithmic
variation in the case of continuous functions. We prescribe a bound on
$\wtv\left(a(\lambda_o)\right)$; for the variation $\tv(v_o, u_o)$
there is not such a bound but, roughly speaking, the larger $\tv(v_o,
u_o)$ is, the smaller must be $\wtv\left(a(\lambda_o)\right)$. An
important point is that we give explicit expressions for these bounds;
then our results are qualitatively different from some of those quoted
above, where a generic smallness is required.

\smallskip

The plan of the paper is the following. The main result is stated in
Section \ref{main}, Theorem \ref{thm:main}.
The Riemann problem is reviewed in Section \ref{prelim} together with
related results; proofs have been given in \cite{AmCo06-Proceed-Lyon}.
The definition of the algorithm is in Section \ref{sec:app_sol}. The
core of the proof are Section \ref{sec:interactions} -- where
interactions are studied in detail -- and Section \ref{sec:Cauchy} --
where we prove the convergence and consistence of the scheme. A
careful analysis is needed due to the presence of large waves.

The paper is completed by two appendices. In the first one we prove
the main result on the weighted total variation. In the second we
study the interaction of two shock waves to the light of Section
\ref{sec:interactions}, namely we look for precise bounds of the
damping coefficient that controls the reflected wave produced in the
interaction; we think that this analysis is interesting by its own.
Good reading!


\section{Main results}\label{main}
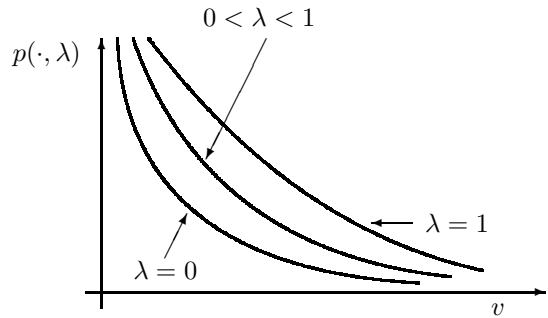
\begin{figure}[htbp]

\begin{picture}(100,120)(0,-10)
\setlength{\unitlength}{1.2pt}


\put(85,0){ \put(0,0){\vector(1,0){140}}
\put(125,-7){\makebox(0,0)[b]{$v$}}
\put(0,0){\line(-1,0){5}}
\put(0,0){\vector(0,1){80}}
\put(-28,75){\makebox(0,0)[l]{$p(\cdot,\lambda)$}}
\put(0,0){\line(0,-1){5}}


\put(0,0){\thicklines{\qbezier(5,80)(8,8)(100,3)}}
\put(10,10){\makebox(0,0)[tl]{$\lambda=0$}}
\put(20,11){\thinlines{\vector(1,2){7}}}
\put(0,0){\thicklines{\qbezier(10,80)(33,13)(110,5)}
\put(52,80){\thinlines{\vector(-1,-2){19}}}
\put(32,90){\makebox(0,0)[tl]{$0<\lambda<1$}}}
\put(0,0){\thicklines{\qbezier(15,80)(60,20)(120,7)}
\put(102,25){\makebox(0,0)[tl]{$\lambda=1$}}}
\put(98,22){\thinlines{\vector(-1,0){13}}}
}

\end{picture}

\caption{\label{fig:pressures}\small{Pressure curves as functions
of $v$.}}

\end{figure}

We consider the system of conservation laws (\ref{eq:system}). The
pressure is given by
\begin{equation}\label{eq:pressure}
p(v,\lambda)= \frac{a^2(\lambda)}v
\end{equation}
where $a$ is a smooth ($\C{1}$) function defined on $[0,1]$ satisfying
for every $\lambda\in[0,1]$
\begin{equation}\label{eq:phi}
a(\lambda)>0, \qquad a'(\lambda)>0\,,
\end{equation}
see Figure \ref{fig:pressures}. For instance $a^2(\lambda) =
k_0+\lambda(k_1-k_0)$ for $0<k_0<k_1$. As a consequence of
(\ref{eq:pressure}) and (\ref{eq:phi}) we have, for every
$(v,\lambda)\in(0,+\infty)\times[0,1]$,
\begin{eqnarray}
p>0,\qquad
&p_v<0,\qquad p_{vv}>0,& \label{eq:pressv}
\\
&p_\lambda>0,\qquad p_{v\lambda}<0\,.& \label{eq:pressl}
\end{eqnarray}
Remark that assumptions (\ref{eq:pressv}) and (\ref{eq:pressl}) are
analogous to those usually made on the pressure in the full
non-isentropic case, \cite{Liu}, the entropy replacing $\lambda$.


We denote $U=(v,u,\lambda)\in\Omega=(0,+\infty) \times \reali \times
[0,1]$ and by $\tilde U=(v,u)$ the projection of $U$ onto the plane
$vu$; the same notation applies to curves. Under assumptions
(\ref{eq:pressure}) and (\ref{eq:phi}) the system (\ref{eq:system}) is
strictly hyperbolic in the whole $\Omega$ with eigenvalues
$e_1=-\sqrt{-p_v(v,\lambda)}$, $e_2=0$, $e_3=\sqrt{-p_v(v,\lambda)}$.
We write $c = \sqrt{-p_v} = a(\lambda)/v$. The eigenvectors associated
to the eigenvalues $e_i$, $i=1,2,3$, are $r_1=(1,c,0)$,
$r_2=(-p_\lambda,0,p_v)$, $r_3=(-1,c,0)$. Because of the third
inequality in (\ref{eq:pressv}) the eigenvalues $e_1$, $e_3$ are
genuinely nonlinear with $\nabla e_i \cdot r_i=p_{vv}/(2c)>0$,
$i=1,3$, while $e_2$ is linearly degenerate. Pairs of Riemann
invariants are $R_1=\{u-a(\lambda)\log v,\lambda\}$, $R_2=\{u,p\}$,
$R_3=\{u+a(\lambda)\log v, \lambda\}$.

\smallskip

We denote by $\tv(f)$ the total variation of a function $f$. In the
case $f:\reali\to(0,+\infty)$ we define the \emph{weighted total
  variation} of $f$ by
\begin{eqnarray}\nonumber
\wtv(f) = 2\sup \sum_{j=1}^n \frac{\left|f(x_j) -
f(x_{j-1})\right|}{f(x_j) + f(x_{j-1})}
\end{eqnarray}
where the supremum is taken over all $n\ge1$ and $(n+1)$-tuples of
points $x_j$ with $x_o<x_1<\ldots<x_n$. This variation is motivated by
the definition (\ref{eq:strengths}) of strength for the waves of the
second family. If $f$ is bounded and bounded away from zero then
clearly
\[
\frac{1}{\sup f}\tv(f)\le \wtv(f)\le \frac{1}{\inf f}\tv(f)\,.
\]

\begin{proposition}\label{prop:wtv}
  Consider $f: \reali \to(0,+\infty)$; then
\begin{equation}\label{eq:smaller}
\frac{\inf f}{\sup f}\tv\left(\log(f)\right)\le \wtv(f)\le
\tv\left(\log(f)\right)\,.
\end{equation}
Moreover, if $f\in C(\reali)$ then $\wtv(f) =
\tv\left(\log(f)\right)$.
\end{proposition}

\smallskip

The proof is deferred to Appendix~\ref{app:WTV}. In
(\ref{eq:smaller}), in the inequality on the right, the strict sign
may occur if $f$ is discontinuous; see Remark~\ref{rem:log-f}.

\smallskip\par

We provide system (\ref{eq:system}) with initial data
\begin{eqnarray}
U(x,0)= U_o(x) = (v_o(x),u_o(x),\lambda_o(x))\label{init-data}
\end{eqnarray}
for $x\in\reali$. Denote $a_o(x) \doteq a\left(\lambda_o(x)\right)$, $p_o(x) \doteq
p\left(v_o(x), \lambda_o(x)\right)$; remark that $\inf a_o(x)\ge a(0)
> 0$. The main result of this paper now follows.

\smallskip

\begin{theorem}\label{thm:main}
  Assume (\ref{eq:pressure}), (\ref{eq:phi}). Consider initial data (\ref{init-data})
  with $v_o(x)\ge \underline{v}>0$ for some constant $\underline{v}$
  and $0\le \lambda_o(x)\le 1$. For every $m > 0$ and a suitable function $k(m) \in (0,1/2)$
  the following holds. If
\begin{eqnarray}
\tv\left(\log(p_o)\right) + \frac{1}{\inf a_o}\tv(u_o) & < &
2\Bigl(1 - 2\wtv(a_o)\Bigr)m \label{hyp1}
\\
\wtv(a_o) & < & k(m) \label{hyp2}
\end{eqnarray}
then the Cauchy problem (\ref{eq:system}), (\ref{init-data}) has a
weak entropic solution $(v,u,\lambda)$ defined for
$t\in\left[0,+\infty\right)$. Moreover the solution is valued in a
compact set of $\Omega$ and there is a constant $C(m)$ such that
for every $t\in\left[0,+\infty\right)$
\begin{eqnarray}\label{eq:tv-est}
\tv \left(v(t,\cdot), u(t,\cdot)\right) & \le & C(m)\,.
\end{eqnarray}
\end{theorem}

The function $k(m)$, whose expression is given in (\ref{eq:kM}),
deserves some comments. The interaction of two waves $\alpha$,
$\alpha'$ of the same family $i=1,3$ produces a wave $\beta$ of the
same family $i$ and a \lq\lq reflected\rq\rq\ wave $\delta$ of the
other family $j$ ($j=1,3$, $j\ne i$). For a suitable definition of the
strengths of the waves we prove that $|\delta| \le d \cdot
\min\{|\alpha|, |\alpha'|\}$ for a damping coefficient $d<1$ depending
on $\alpha$ and $\alpha'$, see Lemma \ref{lem:shock-riflesso}. The
function $k$ above depends essentially on the supremum of such
coefficients $d$; we prove that $k(0)=1/2$ and that $k(m)$ decreases
to $0$ as $m\to+\infty$. In particular then $\wtv(a_o)<1/2$. The
assumptions (\ref{hyp1}), (\ref{hyp2}) read as analogous to those in
\cite{NishidaSmoller}: the larger is $m$, the smaller is $k(m)$, and
vice-versa. The occurring of a possible blow-up when the bound on
$\wtv(a_o)$ does not hold is an interesting open problem.

The variation of $\lambda_o$ appears both in condition
(\ref{hyp1}), because of $p_o$, and in (\ref{hyp2}). Using the
definition of the pressure, we can replace (\ref{hyp1}) by the
slightly stronger condition $\tv \log(v_o) + 2\tv \log(a_o) +
\frac{1}{\inf a_o}\tv(u_o)  \le  2\bigl(1 - 2\wtv(a_o)\bigr)m $ or
even
\[
\tv \log(v_o) + \frac{1}{a(0)}\tv(u_o) \le
2m\left(1-\frac{2m+1}{m}\tv(\log(a_o))\right)
\]
by making use of (\ref{eq:smaller}). In particular if $\lambda_o$
is constant we recover the famous result by Nishida
\cite{Nishida68}.

\smallskip

Clearly $\lambda(t,x)=\lambda_o(x)$ for any $t$ because of the third
equation in (\ref{eq:system}); this is why only $v$ and $u$ appear in
the estimate (\ref{eq:tv-est}). In other words system
(\ref{eq:system}) can be rewritten as a $p$-system of two conservation
laws with flux depending on $x$, namely for the pressure law $p =
p\left(v,\lambda_o(x)\right)=a^2\left(\lambda_o(x)\right)/v$.

\smallskip

The proof of Theorem \ref{thm:main} makes use of a wave-front tracking
scheme where we exploit the special structure of system
(\ref{eq:system}) by differentiating the treatment of $1$ and $3$
waves from that of $2$ waves. Our algorithm is a natural extension of
that in \cite{AmadoriGuerra01}, where the system for $\lambda_o$
constant is studied, in presence of a relaxation term.

Here we consider a linear functional as in \cite{AmadoriGuerra01}
that accounts for the strengths of all $1$ and $3$ waves, with a
weight $\xi>1$ assigned to shock waves; a crucial point in the
proof is the choice of $\xi$ as a function of $m$. This functional
differs with that in \cite{Nishida68}, \cite{Asakura2}, where
$\xi$ is missing and only the variation of shocks is taken into
account. Moreover, motivated again by \cite{AmadoriGuerra01}, we
do not introduce a simplified Riemann solver for interactions
between $1$ and $3$ waves but only for interactions involving the
$2$-contact discontinuities. The interaction potential considers
then uniquely interactions of $2$ waves with $1$ or $3$ waves
approaching to it.

\smallskip

System (\ref{eq:system}) can be written in in Eulerian coordinates.
Denoting $\rho=1/v$ the density, the pressure law becomes
$p=a^2(\lambda)\rho$ and (\ref{eq:system}) turns into
\begin{equation}\label{eq:system-euler}
\left\{
\begin{array}{ll}
\rho_t +(\rho u)_x &= 0\\
(\rho u)_t + \left(\rho u^2+p(\rho,\lambda)\right)_x &= 0
\\
(\rho\lambda)_t +(\rho \lambda u)_x &= 0\,.
\end{array}
\right.
\end{equation}
A global existence result of weak solutions for
(\ref{eq:system-euler}) holds by Theorem \ref{thm:main} because of
\cite{Wagner}.


\section{Preliminaries}\label{prelim}

\subsection{Comparison with other models}\label{sec:comparison}
In \cite{Benzoni-Gavage} many models for diphasic flows are proposed
and studied. In a simple case (no source terms, the fluid in either a
dispersed or separated configuration) and keeping notations as in
\cite[page 35]{Benzoni-Gavage}, they can be written as
\begin{equation}\label{eq:sylvie}
\left\{
\begin{array}{ll}
(\rho_l R_l)_t+(\rho_l R_l u_l)_x & =  0
\\
(\rho_g R_g)_t+(\rho_g R_g u_g)_x & =  0
\\
(\rho_l R_l u_l + \rho_g R_g u_g)_t+(\rho_l R_l u_l^2 + \rho_g R_g
u_g^2+ p)_x &= 0\,.
\end{array}
\right.
\end{equation}
Here the indexes $l$ and $g$ stand for \emph{liquid} and \emph{gas}.
Therefore $\rho_l$, $R_l$, $u_l$ are the liquid density, phase
fraction, velocity, and analogously for the gas; clearly $R_l + R_g =
1$. The pressure law is $p = a^2\rho_g$, for $a>0$ a constant.
Equations (\ref{eq:sylvie}) state the conservation of mass of either
phases and the total momentum.

A case studied in \cite[page 44]{Benzoni-Gavage} is when $u_l =u_g$
and $\rho_l$ is constant, say equal to $1$. The unknown variables are
then $R_l$, $u$, $\rho_g$, and it is assumed $0<R_l<1$ and $\rho_g>0$;
as a consequence $0<R_g<1$. Under these conditions, and writing still
$\rho_l$ instead of $1$ for clarity, we define the concentration $c
=\frac{\rho_g R_g}{\rho_l R_l}>0$ and deduce the pressure law $p = a^2
c \frac{R_l}{1-R_l}$. We obtain exactly the model of \cite{Peng}:
\begin{equation}\label{eq:peng-euler}
\quad \left\{
\begin{array}{ll}
(R_l)_t +(R_l u)_x &= 0\\
(R_l c)_t + \left(R_l c u\right)_x &= 0
\\
(R_l(1+c)u)_t +\left(R_l (1+c)u^2+p\right)_x &= 0\,.
\end{array}
\right.
\end{equation}
This system is strictly hyperbolic for $c>0$. Remark that the three
eigenvalues of (\ref{eq:peng-euler}) coincide with $u$ at $c=0$ and if
$c$ vanishes identically then (\ref{eq:peng-euler}) reduces to the
pressureless gasdynamics system. System (\ref{eq:peng-euler}) is
analogous to (\ref{eq:system-euler}) but the pressure laws are
different. In fact the variables $\rho$ and $\lambda$ of
(\ref{eq:system-euler}) write $\rho=\rho_l R_l + \rho_g R_g$ and
$\lambda = \frac{\rho_g R_g}{\rho_l R_l + \rho_g R_g}=\frac{c}{1+c}$
and then $R_l=(1-\lambda)\rho$, $\rho_g =
\frac{\lambda}{\rho^{-1}-(1-\lambda)}$. If we sum up the two first
equations in (\ref{eq:peng-euler}) we find the first equation in
(\ref{eq:system-euler}); the third (resp. second) equation in
(\ref{eq:peng-euler}) becomes the second (third) equation in
(\ref{eq:system-euler}). The choice $p=a^2\rho_g$ for the pressure in
(\ref{eq:sylvie}) gives $p = a^2
\frac{\lambda}{\frac{1}{\rho}-(1-\lambda)}$.

Notice that the pressure vanishes in presence of a pure liquid phase
and this is the main difference with (\ref{eq:pressure}),
(\ref{eq:phi}).

We compare now (\ref{eq:system-euler}) and (\ref{eq:peng-euler}) in
Lagrangian coordinates. Consider for (\ref{eq:peng-euler}) the change
of coordinates $y= R_l dx-R_l u dt$ based on the streamlines of the
liquid particles (because $R_l=\rho_l R_l$), \cite{Peng}. Denote
$w=\frac{1}{R_l}-1= \frac{R_g}{R_l}=\frac {c}{\rho_g}$. Then for
$p=\frac{a^2 c}{w}$ system (\ref{eq:peng-euler}) turns into
\begin{equation}\label{eq:peng-lagrange}
\left\{
\begin{array}{ll}
w_t -u_x &= 0\\
\left((1+c)u\right)_t + p_x &= 0
\\
c_t & = 0\,.
\end{array}
\right.
\end{equation}
It is more interesting however to consider for system
(\ref{eq:peng-euler}) the change $y = (1+c)R_l dx-(1+c)R_l u dt =\rho
dx-\rho u dt$ into Lagrangian coordinates based on the streamlines of
the \emph{full} density $\rho$. Let $w$ be as above and
$v=\frac{w}{1+c}=\frac{R_g}{\rho}$. Then system (\ref{eq:peng-euler})
becomes system (\ref{eq:system}) with $a^2(\lambda)=a^2\frac{c}{1+c} =
a^2\lambda$. As a consequence the pressure law
$p(v,\lambda)=a^2(\lambda)/v$ does not satisfy (\ref{eq:phi}). This
difficulty can be overcome as follows. Fix any $0 < a_1 < a_2 < a$ and
consider for $c_i=\frac{a_i^2}{a^2-a_i^2}$ the invariant domain
$\{(R_l,u,c)\colon 0< c_1 \le c \le c_2\}$, \cite{Peng}. In this
domain $0<b_1\le\lambda\le b_2<1$, for $b_i=a_i^2/a^2$. If we denote
$\mu=\frac{\lambda-b_1}{b_2-b_1}$ then the function
$b(\mu)=a(\lambda)=a(b_1+(b_2-b_1)\mu)$ makes the pressure law
$p(v,\lambda)=b^2(\mu)/v$, with $\mu\in[0,1]$, satisfy both conditions
in (\ref{eq:phi}).


\subsection{Wave curves and the Riemann problem}\label{sec:waves}
In this section we recall some results about the wave curves for
system (\ref{eq:system}) and the solution to the Riemann problem; see
\cite{AmCo06-Proceed-Lyon} for more details.

The 
shock-rarefaction curves through the point $U_o=(v_o,u_o,\lambda_o)$
for (\ref{eq:system}) are
\begin{eqnarray}
&&\Phi_i(v,U_o)=\left(v,\phi_i(v,U_o),\lambda_o\right)\,,\qquad i=1,3\label{eq:lax13}\\[1mm]
&&\quad\phi_1(v,U_o)  =  \left\{
\begin{array}{lll}
u_o+a(\lambda_o)\cdot(v-v_o)/ {\sqrt{vv_o}}   
& v<v_o& \hbox{shock}
\\[1mm]
u_o+a(\lambda_o)\log({v}/{v_o}) & v>v_o\,, & \hbox{rarefaction,}
\end{array}
\right.
\nonumber
\\[2mm]
&&\quad
\phi_3(v,U_o)  =  \left\{
\begin{array}{lll}
u_o-a(\lambda_o)\log({v}/{v_o})
& v<v_o& \hbox{rarefaction}
\\[2mm]
u_o-a(\lambda_o)\cdot(v-v_o)/ {\sqrt{vv_o}}
& v>v_o\,,& \hbox{shock,}
\end{array}
\right.
\nonumber
\\[1mm]
&&\Phi_2(\lambda,U_o) =
\left(v_o\ds\frac{a^2(\lambda)}{a^2(\lambda_o)},
u_o,\lambda\right)\,,\qquad \lambda\in[0,1] \quad \hbox{contact
discontinuity.} \label{eq:lax2}
\end{eqnarray}

The curves $\Phi_1$, $\Phi_2$ and $\Phi_3$ are \emph{plane} curves:
$\Phi_1$ and $\Phi_3$ lie on the plane $\lambda=\lambda_o$ while
$\Phi_2$ on $u=u_o$.
%
%

\begin{definition}[Wave strengths]
  Under the notations (\ref{eq:lax13}),(\ref{eq:lax2}) we define the
  strength $\eps_i$ of a $i$-wave as
\begin{equation}\label{eq:strengths}
\eps_1=\frac{1}{2}\log\left(\frac{v}{v_o}\right), \quad
\eps_2=2\,\frac{a(\lambda)-a(\lambda_o)}{a(\lambda)+a(\lambda_o)},
\quad \eps_3=\frac{1}{2}\log\left(\frac{v_o}{v}\right)\,.
\end{equation}
\end{definition}

\noindent According to this definition, rarefaction waves have
positive strengths and shock waves have negative strengths. Given the
initial datum $\lambda_o=\lambda_o(x)$,
denote
\begin{equation}\label{eq:a*}
a^* \doteq \sup_{x\in\reali}
a\bigl(\lambda_o(x)\bigr)\,,
\quad
a_* \doteq \inf_{x\in\reali} a\bigl(\lambda_o(x)\bigr)\,,
\quad
[a]_* \doteq \frac{a^* - a_*}{a^* + a_*}\,.
\end{equation}
Then $[a]_* \le \frac{a(1) - a(0)}{a(1) + a(0)} < 1$ and $|\eps_2|\le
2[a]_* < 2$. It is useful to define also the function, see
\cite{Peng},
\begin{equation}
\label{h}
h(\eps)=\left\{
\begin{array}{ll}
\eps&\hbox{ if }\eps \ge 0\,,
\\
\sinh \eps& \hbox{ if }\eps < 0\,.
\end{array}
\right.
\end{equation}
Then we have for $i=1,3$
\begin{equation}\label{eq:phis}
\phi_i(v,U_o)= u_o + a(\lambda_o) \cdot2 h(\eps_i)\,.
\end{equation}


At last we consider the Riemann problem. This is the initial-value
problem for (\ref{eq:system}) under the piecewise constant initial
condition
\begin{equation}\label{eq:incond}
(v,u,\lambda)(0,x)=\left\{
\begin{array}{ll}
(v_\ell,u_\ell,\lambda_\ell)=U_\ell & \hbox{ if }x<0
\\
(v_r,u_r,\lambda_r)=U_r & \hbox{ if }x>0
\end{array}
\right.
\end{equation}
for $U_\ell$ and $U_r$ in $\Omega$. We denote $a_r=a(\lambda_r)$,
$p_r= a^2_r/v_r$, and similarly $a_\ell$, $p_\ell$.

\smallskip

\begin{proposition}
\label{prop:RP}
Fix any pair of states $U_\ell$, $U_r$ in $\Omega$; then the Riemann
problem (\ref{eq:system}), (\ref{eq:incond}) has a unique
$\Omega$-valued solution in the class of solutions consisting of
simple Lax waves. If $\eps_i$ is the strength of the $i$-wave,
$i=1,2,3$, then
\[
\eps_3-\eps_1  =
\frac{1}{2}\log\left(\frac{p_r}{p_\ell}\right), \qquad
2\bigl(a_\ell h(\eps_1) +a_r h(\eps_3)\bigr)   =
u_r-u_\ell\,.
\]
Moreover, let $\underline{v}>0$ be a fixed number. There exists a
constant $C_1>0$ depending on $\underline{v}$ and $a(\lambda)$
such that if ${\tilde U}_l$, ${\tilde U}_r\in\Omega$ and
$v>\underline{v}$, then
\begin{eqnarray}
|\eps_1| + |\eps_2| + |\eps_3| \leq C_1 | U_\ell
 - U_r|\,.\label{stima-ampiezze-RP}
\end{eqnarray}
\end{proposition}

For the proof, see \cite{AmCo06-Proceed-Lyon}. One can easily find
that
\begin{eqnarray}\label{stima-onde-PdR}
&&\qquad\left| \eps_1 \right| + \left|\eps_3\right|
~ \leq ~
\frac{1}{2}|\log(p_r) -\log(p_\ell)|
 + \frac 1 {2
\min\{a_\ell,a_r\}}\left|u_r-u_\ell \right|\\\nonumber
&&\qquad\qquad \leq~\frac 1 2 |\log(v_r) -\log(v_\ell)| +  |\log(a_r) -\log(a_\ell)|
 + \frac 1 {2
\min\{a_\ell,a_r\}}\left|u_r-u_\ell \right|\,.
\end{eqnarray}
We remark that for any Riemann data
$(v_\ell,u_\ell,\lambda_\ell)$, $(v_r,u_r,\lambda_r)$, the
$\lambda$ component of the solution takes value $\lambda_\ell$ for
$x<0$ and $\lambda_r$ for $x>0$. The fact that the interfaces
between different phases are connected by a stationary wave can be
interpreted then as a \lq\lq kinetic condition\rq\rq,
\cite{Abeyaratne-Knowles}, analogous to Maxwell's rule.


\section{The approximate solution}\label{sec:app_sol}
In this section we define a wave-front tracking scheme
\cite{Bressanbook} to build up piecewise constant approximate
solutions to (\ref{eq:system}). More precisely we follow the algorithm
introduced in \cite{AmadoriGuerra01}.

First, we approximate the initial data. For any $\nu\in\naturali$ we
take a sequence $(v^\nu_o,u^\nu_o,\lambda^\nu_o)$ of piecewise
constant functions with a finite number of jumps such that

\smallskip
\begin{romannum}
\item $\tv p^\nu_o\leq \tv p_o$, $\tv u^\nu_o\leq \tv u_o$, $\wtv
  a(\lambda^\nu_o) \leq \wtv a(\lambda_o)$, $\inf
  a_o^\nu \ge \inf a_o$;

\item $\lim_{x\to-\infty} (v^\nu_o,u^\nu_o,\lambda^\nu_o)(x)
  =\lim_{x\to-\infty} (v_o,u_o,\lambda_o)(x)$;

\item $\|(v^\nu_o,u^\nu_o,\lambda^\nu_o) -
(v_o,u_o,\lambda_o)\|_{\L1}\leq \frac 1 \nu$
\end{romannum}

\smallskip\par\noindent where $p^\nu_o=a^2(\lambda^\nu_o)/v^\nu_o$.
Second, we define the approximate Riemann solver. We introduce
positive parameters $\eta=\eta_\nu$, $\rho=\rho_\nu$;
they control respectively the size of rarefactions and the threshold
when a simplified Riemann solver is used. Define also a parameter
$\hat s>0$ strictly larger than all possible speeds of wave-fronts of
both families $1$ and $3$. These parameters will be determined at the
end of Section~\ref{sec:Cauchy}.

\smallskip

\begin{itemize}
\item At time $t=0$ we solve the Riemann problems at each point of
  jump of $(v^\nu_o, u^\nu_o, \lambda^\nu_o)(0+, \cdot\,)$ as follows:
  shocks are not modified while rarefactions are approximated by fans
  of waves, each of them having size less than $\eta$. More precisely,
  a rarefaction of size $\eps$ is approximated by $N=[\eps/\eta]+1$
  waves whose size is $\eps/N<\eta$; we set their speeds to be equal
  to the characteristic speed of the state at the right.

  Then $(v,u,\lambda)(t,\cdot)$ is defined until some wave fronts
  interact; by slightly changing the speed of some waves
  \cite{Bressanbook} we can assume that only \emph{two} fronts
  interact at a time.

\item When two wave fronts of families either $1$ or $3$ interact we
  solve the Riemann problem at the interaction point. If one of the
  incoming waves is a rarefaction, after the interaction it is
  prolonged (if it still exists) as a single discontinuity with speed
  equal to the characteristic speed of the state at the right.  If a
  new rarefaction is generated, we employ the Riemann solver described
  before and divide it into a fan of waves having size less than
  $\eta$.

\item When a wave front either of family either $1$ or $3$ interacts
  with a 2-wave we proceed as follows. Let $\delta_2$ be the size of
  the 2-wave and $\delta$ the size of the other wave.

\begin{itemize}
\item If $|\delta_2\delta|\geq \rho$ we solve the Riemann problem as
  above, that is with the \emph{accurate Riemann solver}.
\item If $|\delta_2\delta|< \rho$ we prolong the 1- or 3- wave with a
  wave of the same family and size. Since the two waves do not
  commute, a \emph{non-physical} front is introduced,
  \cite{Bressanbook}, with fixed speed $\hat s>0$. The size of a
  non-physical wave is set to be $|u_r - u_\ell|$, where $u_\ell$,
  $u_r$ are the $u$ components of the left and right states of the
  wave. We call this solver the \emph{simplified Riemann solver}.
\end{itemize}

\item When a non-physical front interacts with a front of family $1$,
  $2$ or $3$ (``physical''), we prolong the solution with a physical
  wave of the same size and a non-physical one, computing the
  intermediate value consequently.

\end{itemize}

\smallskip

We refer for the last two items to Proposition~\ref{lem:np} below.
Remark that two non-physical front cannot interact since they have the
same constant speed $\hat s$. We denote by ${\cal N\!P}$ the set
of non-physical waves.


\section{Interactions}\label{sec:interactions}
Fix the index $\nu$ introduced in the previous section. We shall
prove in Subsection \ref{subs:number_int} that the algorithm
described above is defined for any $t>0$ and provides for any
initial data $(v_o^\nu,u_o^\nu,\lambda_o^\nu)$ a piecewise
constant approximate solution $(v^\nu,u^\nu,\lambda^\nu) =
(v,u,\lambda)$, where we dropped for simplicity the index $\nu$.
Here we study the interaction of waves.

For $\knp>0$ and $t>0$ we define the functional $L$ and the
interaction potential $Q$, both referred to $(v,u,\lambda)(t,\cdot)$,
by
\begin{eqnarray}
L(t) & = & \sum_{i=1,3}|\gamma_i| + \knp L_{np}\,,\qquad\qquad
L_{np} = \sum_{\gamma \in {\cal N\!P}} |\gamma|
\nonumber\\
Q(t) &=& \sum_{\gamma_3 \mbox{\sz \ at the left of } \delta_2}
|\gamma_3| |\delta_2| + \sum_{\gamma_1 \mbox{\sz\ at the right of
} \delta_2} |\delta_2| |\gamma_1|\,.\label{Q}
\end{eqnarray}
Remark that $L$ takes only into account the strengths of both $1$
and $3$ waves and that of non-physical waves. For contact
discontinuities we define
\[
L_{cd} = \sum|\gamma_2| = \wtv a(\lambda^\nu_o)\,.
\]
Finally, for $\xi\ge1$ and $K\ge0$ we introduce
\begin{eqnarray}\label{L-xi}
L_{\xi} & = & L_{{\rm rarefactions}}+\xi L_{{\rm shocks}} +
\knp L_{np}\,,
\\
\label{F}
F & = & L_{\xi}+ K Q\,.
\end{eqnarray}
For simplicity we omitted to note the dependence on $\knp$ in the
functional $L_\xi$ and on $\knp$, $\xi$, $K$ in $F$; the choice of
$\knp$ shall depend on that of $K$, see Proposition~\ref{lem:np}.

Observe that, if $\lambda_o$ is constant, then $Q=0$ and $F =
L_{\xi}$, whose variation was analyzed in Lemma~3.2 of
\cite{AmadoriGuerra01}. Hence we will assume from now on that
\begin{eqnarray}\label{eq:Ao>0}
A_o \doteq \wtv(a_o) & > &0\,.
\end{eqnarray}
By assumption (i) in Section~\ref{sec:app_sol}, one has $L_{cd}\le
A_o$.

\begin{figure}[htbp]
\begin{picture}(100,80)(-120,0)
\setlength{\unitlength}{1pt}

\put(50,0){

\put(-40,0){\line(1,2){20}} \put(-40,0){\line(0,1){40}}
\put(-40,0){\line(-1,1){40}}
\put(-62,30){\makebox(0,0){$\alpha_1$}}
\put(-47,32){\makebox(0,0){$\alpha_2$}}
\put(-18,30){\makebox(0,0){$\alpha_3$}}

\put(60,0){\line(1,2){20}}
\put(60,0){\line(0,1){40}}
\put(60,0){\line(-3,5){24}}
\put(37,30){\makebox(0,0){$\beta_1$}}
\put(55,33){\makebox(0,0){$\beta_2$}}
\put(83,30){\makebox(0,0){$\beta_3$}}

\put(10,40){\line(2,3){26}} \put(10,40){\line(0,1){40}}
\put(10,40){\line(-3,5){24}} \put(-13,70){\makebox(0,0){$\eps_1$}}
\put(17,73){\makebox(0,0){$\eps_2$}}
\put(38,70){\makebox(0,0){$\eps_3$}}

\put(-80,10){\makebox(0,0){$\ell$}}
\put(10,10){\makebox(0,0){$m$}} \put(100,10){\makebox(0,0){$r$}}

}

\end{picture}

\caption{\label{fig:Glimmone}{A general interaction pattern.}}

\end{figure}
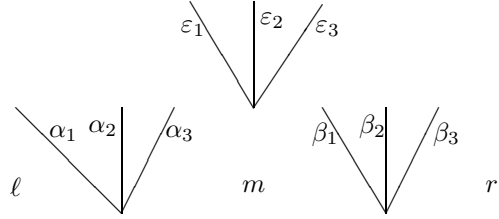

In the following sections we analyze in detail the different types of
interactions. Recalling the definition of $h$, (\ref{h}), and with the
notation of Figure \ref{fig:Glimmone}, we introduce the following
identities, see (3.1), (3.2) in \cite{Peng}:
%
\begin{eqnarray}
\eps_3 -  \eps_1 & = &  \alpha_3 + \beta_3 -\alpha_1 - \beta_1
\label{tre-uno}
\\
a_\ell h(\eps_1) + a_r h(\eps_3) & = &  a_\ell h(\alpha_1) + a_m
h(\alpha_3) + a_m h(\beta_1) + a_r h(\beta_3)\,. \label{tre-due}
\end{eqnarray}


\noindent Formula (\ref{tre-uno}) does not depend on $\lambda$ and
follows easily by equating the specific volumes $v$ before and after
the interaction time. By equating the velocities $u$ we obtain
(\ref{tre-due}). These properties are a consequence of the definition
(\ref{eq:strengths}) of the strengths for 1- and 3- waves and of
(\ref{eq:phis}).


\subsection{Interactions with a $2$-wave}
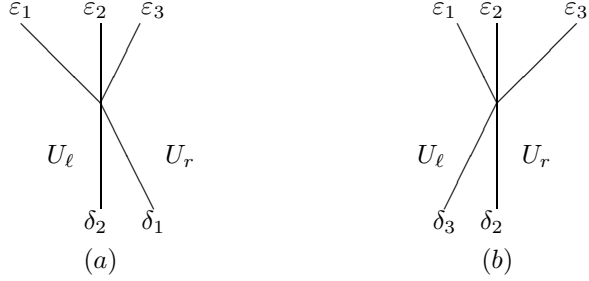
\begin{figure}[htbp]
\begin{picture}(100,100)(-80,-15)
\setlength{\unitlength}{1pt}


\put(30,0){

\put(0,0){\line(0,1){40}} \put(-2,-5){\makebox(0,0){$\delta_2$}}
\put(20,0){\line(-1,2){20}} \put(20,-5){\makebox(0,0){$\delta_1$}}
\put(0,40){\line(0,1){30}} \put(-2,75){\makebox(0,0){$\eps_2$}}
\put(0,40){\line(1,2){15}} \put(20,75){\makebox(0,0){$\eps_3$}}
\put(0,40){\line(-1,1){30}} \put(-30,75){\makebox(0,0){$\eps_1$}}
\put(-15,20){\makebox(0,0){$U_\ell$}}
\put(30,20){\makebox(0,0){$U_r$}}

\put(0,-20){\makebox(0,0){$(a)$}} }

\put(180,0){

\put(0,0){\line(0,1){40}} \put(-2,-5){\makebox(0,0){$\delta_2$}}

\put(-20,0){\line(1,2){20}}
\put(-20,-5){\makebox(0,0){$\delta_3$}}

\put(0,40){\line(0,1){30}} \put(-2,75){\makebox(0,0){$\eps_2$}}

\put(0,40){\line(1,1){30}} \put(30,75){\makebox(0,0){$\eps_3$}}

\put(0,40){\line(-1,2){15}} \put(-20,75){\makebox(0,0){$\eps_1$}}

\put(-25,20){\makebox(0,0){$U_\ell$}}
\put(15,20){\makebox(0,0){$U_r$}}

\put(0,-20){\makebox(0,0){$(b)$}} }

\end{picture}

\caption{\label{fig:inter}{Interactions. $(a)$: from the right;
$(b)$: from the left.}}
\end{figure}

We consider first the interactions of $1$ or $3$ waves with a $2$
wave, see Figure \ref{fig:inter}.

\smallskip

\begin{proposition}[\cite{AmCo06-Proceed-Lyon}]\label{interazzioni}
  Denote by $\lambda_\ell$, $\lambda_r$ the side states of a $2$-wave.
  The interactions of $1$- or $3$-waves with the $2$-wave give rise to
  the following pattern of solutions:
\[
\begin{tabular}{c|cc}
interaction & \multicolumn{2}{c}{outcome}
\\
\cline{2-3} & $\lambda_\ell<\lambda_r$ & $\lambda_\ell>\lambda_r$
\\
\hline $2\times 1R$ & $1R+2+3R$ & $1R+2+3S$
\\
$2\times 1S$ & $1S+2+3S$ & $1S+2+3R$
\\
$3R\times 2$ & $1S+2+3R$ & $1R+2+3R$
\\
$3S\times 2$ & $1R+2+3S$ & $1S+2+3S$.
\end{tabular}
\]
\end{proposition}

The next lemma is concerned instead with the \emph{strengths} of
waves involved in the interaction above. The inequalities
(\ref{eq:stima-interazione-semplice}) improve the inequality (3.3)
in \cite{Peng} in the special case of two interacting wave fronts,
one of them being of the second family. More precisely under the
notations of \cite{Peng} we find a term $1/(a_r+a_\ell)$ instead
of $1/\min\{a_r, a_\ell\}$. The proof differs from Peng's. Our
estimates are sharp: in some cases
(\ref{eq:stima-interazione-semplice}) reduces to an identity.


\begin{lemma}[\cite{AmCo06-Proceed-Lyon}]\label{lem:interazioni}
Assume that a $1$ wave of strength $\delta_1$ or a $3$ wave of
strength $\delta_3$ interacts with a $2$ wave of strength
$\delta_2=2(a_r - a_\ell)/(a_r+a_\ell)$. Then, the strengths
$\eps_i$ of the outgoing waves satisfy: $\eps_2=\delta_2$ and
\begin{eqnarray}
|\eps_i - \delta_i| = |\eps_j| & \le &  \displaystyle
\frac{1}{2}\,|\delta_2|\cdot |\delta_i| ~\leq~ [a]_*  |\delta_i|\label{eq:33Pengnew}
\end{eqnarray}
for $i,j=1,3$, $i\ne j$. Moreover,
\begin{eqnarray}
|\eps_1|+|\eps_3| & \leq & \left\{
\begin{array}{ll}
|\delta_1| + |\delta_1|\piu{\delta_2}&\qquad\hbox{ if $1$
interacts}
\\
|\delta_3| + |\delta_3|\meno{\delta_2}&\qquad\hbox{ if $3$
interacts\,.}
\end{array}
\right.
\label{eq:stima-interazione-semplice}
\end{eqnarray}
\end{lemma}

Here $[x]_+=\max\{x,0\}$, $[x]_-=\max\{-x,0\}$, $x\in\reali$.
Remark that the colliding $1$ or $3$ wave does not change sign
across the interaction. Moreover the functional $L$ increases if
and only if the incoming and the reflected waves are of the same
type; this happens when the colliding wave is moving toward a more
liquid phase.

Now we prove that $F$ is decreasing for suitable $K$ when an
interaction with a $2$-wave occurs. The potential $Q$ is needed to
balance the possible increase of $L_\xi$.

\smallskip

\begin{proposition}\label{Delta-F-2wave} Assume  $A_o < 2$ and
  consider an interaction of a $1$ or $3$ wave with a $2$ wave, with
  the notation of Lemma~\ref{lem:interazioni}. Then $\Delta Q<0$. If
  moreover
\begin{eqnarray}\label{xi_e_K}
\xi\ge 1\qquad {and}\qquad K > \frac{2\xi}{2-A_o}
\end{eqnarray}
then
\begin{eqnarray}\label{Kappa-mezzi}
\xi|\eps_j| = \xi\bigl||\eps_i| - |\delta_i|\bigr| < \frac K2 |\Delta Q|
\end{eqnarray}
and hence $\Delta F < 0$.
\end{proposition}

\smallskip

\begin{proof}
  We consider the interaction of a $3$-wave with a $2$-wave, as in the
  proof of Lemma~\ref{lem:interazioni},
  see~\cite{AmCo06-Proceed-Lyon}; the symmetric case follows in an
  analogous way. We use the notation as in Figure~\ref{fig:inter}(a).
  We define $L_{cd}^*=L_{cd}^-+L_{cd}^+$, $L_{cd}^\pm$ meaning {\em right or left}
  of the $2$ wave under consideration.

  By assumption, one has
\begin{eqnarray}\label{bound-su-wtv}
L_{cd}=L_{cd}^-+L_{cd}^+ + |\delta_2| = L_{cd}^* + |\delta_2|\leq
A_o < 2\,.
\end{eqnarray}
Recall that $\eps_1-\delta_1 = \eps_3$ and
\begin{eqnarray}
|\eps_1|-|\delta_1| & = & |\eps_3|\,,\quad\quad  \hbox{ if }
\delta_2>0\,,\label{eq:12+}
\\
|\eps_1|-|\delta_1| & = & - |\eps_3|\,,\quad\ \hbox{ if }
\delta_2<0\,,\label{eq:12-}
\end{eqnarray}
so that in particular $|\eps_3| = \bigl||\eps_1| -
|\delta_1|\bigr|$. An estimate for $\Delta Q$ follows at once
because of (\ref{bound-su-wtv}):
\begin{eqnarray}
\Delta Q & = & - |\delta_2\delta_1| +
\left(|\eps_1|-|\delta_1|\right)L_{cd}^- + |\eps_3|L_{cd}^+
 \, \leq\,  \frac{1}{2}|\delta_2\delta_1|(L_{cd}^*-2)
\nonumber
\\
& \leq & \frac{1}{2}|\delta_2\delta_1|(A_o-2)<0\,.
\label{eq:DeltaQ}
\end{eqnarray}
Hence, using (\ref{eq:33Pengnew}), we get
\begin{eqnarray}
\xi|\eps_3| + \frac K 2 \Delta Q \leq
\frac{1}{2}|\delta_2\delta_1|\left\{\xi + \frac K 2  (A_o-2)
\right\}<0 \label{Kappa-mezzib}
\end{eqnarray}
because of (\ref{xi_e_K}); this proves (\ref{Kappa-mezzi}). Finally,
by using (\ref{Kappa-mezzib}) we get
\begin{eqnarray}\label{diseq:F-xi-K}
\Delta F = \Delta L_\xi + K  \Delta Q \leq \xi|\eps_3| + \xi
|\eps_1 - \delta_1| + K \Delta Q <0\,.
\end{eqnarray}
\end{proof}


\subsection{Interactions between $1$ and $3$ waves}
\begin{figure}[htbp]
\begin{picture}(100,80)(-80,-15)
\setlength{\unitlength}{0.8pt}

\put(30,0){
\put(0,40){\line(-2,-3){30}}\put(-40,0){\makebox(0,0){$\delta_3$}}
\put(20,0){\line(-1,2){20}} \put(30,0){\makebox(0,0){$\delta_1$}}
\put(0,40){\line(1,2){15}} \put(20,75){\makebox(0,0){$\eps_3$}}
\put(0,40){\line(-1,1){30}} \put(-30,75){\makebox(0,0){$\eps_1$}}
\put(0,-20){\makebox(0,0){$(a)$}}
}

\put(240,0){ \put(0,40){\line(-1,-1){45}}
\put(-55,0){\makebox(0,0){$\alpha_3$}}
\put(0,40){\line(-2,-3){30}}\put(-10,0){\makebox(0,0){$\beta_3$}}
\put(0,40){\line(1,2){15}} \put(20,75){\makebox(0,0){$\eps_3$}}
\put(0,40){\line(-1,1){30}} \put(-30,75){\makebox(0,0){$\eps_1$}}
\put(0,-20){\makebox(0,0){$(b)$}}
}

\end{picture}

\caption{\label{fig:inter3133}{Interactions of $1$ and $3$
waves.}}
\end{figure}
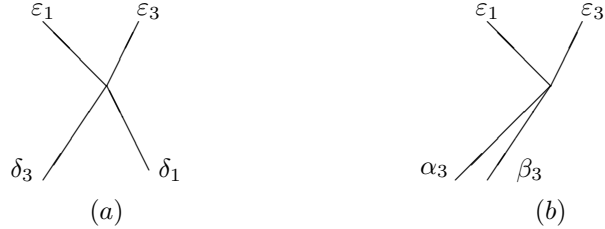

Here we analyze the possible interactions between $1$- and $3$-waves.
Two situations may occur, see Figure \ref{fig:inter3133}: either the
waves belong to different families or they both belong to the same
family. In this last case, at least one of the waves must be a shock.


\begin{lemma}[Different families interacting]\label{diff-family}
  If a wave of the third family interacts with a wave of the first
  family, they cross each other without changing their strength.
\end{lemma}

\begin{proof} See also Lemma~3.1 in \cite{AmadoriGuerra01}.
  Using notation as in Figure~\ref{fig:inter3133}(a) we have
  $\eps_3-\eps_1 = \delta_3-\delta_1$ and $h(\eps_1)+ h(\eps_3) =
  h(\delta_1) + h(\delta_3)$. The uniqueness of solutions to the
  Riemann problem implies $\eps_1=\delta_1$, $\eps_3=\delta_3$.

  Remark that here $\Delta L_{\xi}=0=\Delta Q$ and then $ \Delta F=0$ for
  all $\xi\ge 1$ and $K$.
\end{proof}

\smallskip

\begin{lemma}[Same family interacting: outcome]\label{same-family}
  Assume that a wave $\alpha_3$ of the third family interacts with a
  wave $\beta_3$ of the third family, giving rise to waves $\eps_1$,
  $\eps_3$. Then
\begin{romannum}
\item $\alpha_3<0$, $\beta_3<0$ $\Rightarrow$ $\eps_1>0$, $\eps_3<0$,
\item $\alpha_3\beta_3<0$ $\Rightarrow$ $\eps_1<0$.
\end{romannum}
An analogous result holds for interacting waves of the first family.
\end{lemma}

\smallskip

\begin{proof} The proof can be done in a geometric way by
  observing the mutual positions of the curves, \cite{Nishida68,
    Smoller}. A simple alternative proof by analytical arguments now
  follows. We have
\begin{eqnarray}
\eps_3-\eps_1 & = & \alpha_3+\beta_3
\label{1-0s}\\
h(\eps_1)+ h(\eps_3) & = & h(\alpha_3) + h(\beta_3)\,.
\label{2-0s}
\end{eqnarray}
In case {\em(i)} these formulas read $\eps_1-\eps_3 =
|\alpha_3|+|\beta_3|>0$ and $-h(\eps_1)- h(\eps_3) = \sinh(|\alpha_3|)
+ \sinh(|\beta_3|)>0$. If it were $\eps_3>0$ then $\eps_1>0$ from the
first equality and $\eps_1<0$ from the second, a contradiction.
Therefore $\eps_3<0$ so that $\eps_1+|\eps_3| = |\alpha_3|+|\beta_3|$
and $ -h(\eps_1)+ \sinh(|\eps_3|) = \sinh(|\alpha_3|) +
\sinh(|\beta_3|)$.  Analogously, if it were $\eps_1<0$, using
elementary inequalities we get $0=\sinh(|\eps_1|)+ \sinh(
|\alpha_3|+|\beta_3| + |\eps_1| ) - \sinh(|\alpha_3|) -
\sinh(|\beta_3|) \geq 2 \sinh(|\eps_1|)$, a contradiction again. Hence
$\eps_1>0$.

In case {\em(ii)} assume $\alpha_3<0$, $\beta_3>0$; the other case is
dealt analogously since (\ref{1-0s}), (\ref{2-0s}) are symmetric in
$\alpha_3$, $\beta_3$. We have $\eps_3-\eps_1 = -|\alpha_3|+|\beta_3|$
and $h(\eps_1)+ h(\eps_3) = -\sinh(|\alpha_3|) + |\beta_3|$, then
$[h(\eps_1)+ \eps_1]+ [h(\eps_3)- \eps_3] =
|\alpha_3|-\sinh(|\alpha_3|)<0$. If $\eps_3>0$, this last equality
becomes $h(\eps_1)+ \eps_1 = |\alpha_3|-\sinh(|\alpha_3|)<0$ that
implies $\eps_1<0$. If $\eps_3<0$, then $h(\eps_1)+ \eps_1 =
[|\alpha_3|-\sinh(|\alpha_3|)] - [|\eps_3|-\sinh(|\eps_3|)]$.  If it
were $\eps_1>0$ it would be $|\alpha_3|<|\eps_3|$, since the map
$x\mapsto x-\sinh x$ is decreasing; but from $|\eps_3|+|\eps_1| =
|\alpha_3|-|\beta_3|$ we would get that $|\eps_3|<|\alpha_3|$, a
contradiction. Hence in all cases one has $\eps_1<0$.
\end{proof}

\smallskip

Now we give sharper estimates for the interaction of waves of the
same family: we prove that the strength of the reflected wave is
bounded by the size of each incoming wave, multiplied by a damping
factor smaller than 1. This property will be crucial in the next
section and it holds also for interactions with a 2-wave, with
damping factor $[a]_*$, see (\ref{eq:33Pengnew}). In the case
below however the coefficient depends on the strengths of the
incoming waves; this happens also when non-physical waves are
generated, see Proposition \ref{lem:np}. We assume that
\begin{equation}\label{rogna}
\begin{array}{l}
\hbox{\em the strength of any interacting $i$-wave is less than
$m$,}
\\
\hbox{\em for some $m>0$ and $i=1,3$.}
\end{array}\end{equation}
In the special case of interaction of waves of the same family
producing two outgoing shocks we give a more precise result in
Appendix~\ref{app:SS}.

\smallskip

\begin{lemma}[Same family interacting]\label{lem:shock-riflesso}
  Consider the interaction of two waves of the same family, of sizes
  $\alpha_i$ and $\beta_i$, $i=1,3$, producing two outgoing waves
  $\eps_1$, $\eps_3$; assume (\ref{rogna}). Then the following holds.
\begin{romannum}
\item There exists a damping coefficient $d=d(m)$, with $0<d<1$, such
  that
\begin{equation}\label{eq:chi_def}
|\eps_j| \le d(m) \cdot \min\{|\alpha_i|,|\beta_i|\}\,,\qquad j\ne
i\,.
\end{equation}

\item If the incoming waves are both shocks, the resulting shock
satisfies $|\eps_i|>\max \{|\alpha_i|,|\beta_i|\}$. If the
incoming waves have different signs, both the amount of shocks and
the amount of rarefactions of the $i^{th}$ family decrease across
the interaction.

In any case
  \begin{eqnarray}
  |\eps_i| & \le & |\alpha_i| + |\beta_i|\,.
  \label{eq:333}
  \end{eqnarray}
\end{romannum}
\end{lemma}

\begin{proof} {\it (i)} To fix the ideas, assume $i=3$. We have
\begin{eqnarray}
\eps_3-\eps_1 & = & \alpha_3+\beta_3
\label{1-0}\\
h(\eps_1)+ h(\eps_3) & = & h(\alpha_3) + h(\beta_3) \label{2-0}
\end{eqnarray}
and then
\begin{eqnarray}\label{eq:identita-eps1-alpha3-beta3}
h(\eps_1)+ h(\eps_1+\alpha_3+\beta_3) & = & h(\alpha_3) + h(\beta_3)\,.
\end{eqnarray}
Remark that the equation (\ref{eq:identita-eps1-alpha3-beta3}) is
symmetric in $\alpha_3$, $\beta_3$; from the implicit function theorem
we find that $\eps_1=\eps_1(\alpha_3,\beta_3)$ is $\C1$.  Using the
notation $\eps_1=\t$, $\alpha_3=a$, $\beta_3=b$, the identity
(\ref{eq:identita-eps1-alpha3-beta3}) rewrites as
\begin{eqnarray}\label{eq:tau-a-b}
h(\t)+ h(\t+a+b) - h(a) - h(b)=0\,,
\end{eqnarray}
with $\t=\t(a,b)$. One verifies that $\t(a,0)=\t(0,b)=0$ and that
$$
\t_a = \frac{h'(a) - h'(\t+a+b)}{h'(\t) + h'(\t+a+b)},\qquad \t_b =
\frac{h'(b) - h'(\t+a+b)}{h'(\t) + h'(\t+a+b)}\,.
$$
As $a\to 0$, one has $\t_a\to (1-h'(b))/ (1+h'(b))$; then
$|\t_a(0,b)|<1$ and it can be bounded by a positive constant less
than 1 that depends on $m$. The same argument works for $\t_b$.

To complete the proof, we show that $|\t| < \min\{|a|,|b|\}$, in
the non-trivial case $a\not =0\not=b$. We argue by contradiction,
using an argument of \cite{Poupaud-Rascle-Vila}. Suppose that
$|\t| \geq |a|$; we can assume $\t>0$, since the case $\t<0$ can
be proved by using the equality (\ref{eq:tau-a-b}) written in
terms of $G(t)=-h(-t)$. Since the function $h$ is increasing we
have, for $\t\geq |a|$,
$$
h(\t)\geq  h(a),\qquad  h(\t+a+b)\geq h(b)
$$
Moreover, one of the two inequalities is strict: if $a<0$ the
first, if $a>0$ the second. Hence we contradict
(\ref{eq:tau-a-b}).

\smallskip\par {\it (ii)} From \cite{AmadoriGuerra01} we already
know that $\Delta L = |\eps_1| +  |\eps_3| - |\alpha_3| -
|\beta_3|\leq 0$\,, hence (\ref{eq:333}).

If the incoming waves are both shocks, then (\ref{1-0}) becomes
\begin{eqnarray}\label{1-SS>RS}
|\eps_3|+|\eps_1|  & = &  |\alpha_3|+|\beta_3|\,.
\end{eqnarray}
From {\it (i)} we have $|\eps_1| < |\alpha_3|$, $|\beta_3|$ and
hence the first part of {\it (ii)}. On the other hand, if
$\alpha_3\beta_3<0$, we have $\eps_1<0$. We can assume
$\alpha_3<0<\beta_3$; hence (\ref{1-0}) becomes $\eps_3 =
|\beta_3| - |\alpha_3| - |\eps_1|$. If $\eps_3>0$, then $|\eps_3|<
|\beta_3|$; if $\eps_3<0$, using {\it (i)} again one finds
$|\eps_3|< |\alpha_3|$.
\end{proof}

\medskip

\begin{figure}[htbp]
  \centering
  \begin{psfrags}
    \psfrag{d(m)}{\small$\ d(m)$}
    \psfrag{m}{\small$\ m$}
    \hspace{-5mm}\includegraphics[width=6cm]{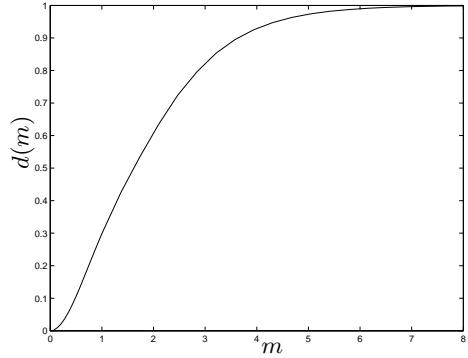}
  \end{psfrags}
  \caption{\label{fig:figure}\small{The coefficient $d(m)$.}}
\end{figure}

\begin{remark}\label{rem:d(m)}\rm The damping coefficient $d(m)$, see Figure \ref{fig:figure},
is given by
\begin{eqnarray*}
d(m) & = & \max_{|a|\le m\atop |b|\le
m}\frac{|\eps(a,b)|}{\min\{|a|, |b|\}}\,,
\end{eqnarray*}
where the function $\eps(a,b)$ satisfies $h(\eps) + h(\eps + a +
b) - h(a) - h(b)=0$, see (\ref{eq:identita-eps1-alpha3-beta3}).
Hence $d(m)$ increases with $m$, and vanishes as $m\to0$ because
quadratic interaction estimates hold for $m$ small.

Moreover, it is asymptotic to $1$ for $m$ large. Indeed, from the
proof of Lemma~\ref{lem:shock-riflesso} we have $\tau_a(0,b) =
\frac{1-h'(b)}{1+h'(b)}$; then $\tau_a(0,b)=0$ if $b>0$ and
$|\tau_a(0,b)|=\frac{\cosh(b)-1}{\cosh(b)+1} \le
\frac{\cosh{m}-1}{\cosh{m}+1}$ if $b<0$. Therefore $|\tau_a(0,b)|
\le \frac{\cosh{m}-1}{\cosh{m}+1}$ for every $b$, and an analogous
estimate holds for $|\tau_b(0,a)|$. Hence $d(m) \ge
\frac{\cosh{m}-1}{\cosh{m}+1} \doteq c(m)$; we refer to Lemma
\ref{lem:second} for the role of this quantity.
\end{remark}

\begin{remark}\label{rem:rar_decreasing}\rm
  When a rarefaction interacts with a 1- or 3-wave, its size does not
  increase. Indeed, the size does not change upon interactions with
  waves of the other family by Lemma \ref{diff-family}; if the
  rarefaction interacts with a shock of the same family we apply
  Lemma~\ref{lem:shock-riflesso}\emph{(ii)}. Remark moreover that by
  Lemmas \ref{diff-family}, \ref{same-family} a rarefaction never
  produces a rarefaction of the other family by interactions with $1$-
  and $3$-waves.
\end{remark}

\smallskip

\begin{remark}\label{rem:out-waves}\rm
  If two waves of the same family interact, the wave belonging to that
  family can be missing, while the \lq\lq reflected\rq\rq\ wave is
  always present. This follows easily from (\ref{1-0}),
  (\ref{2-0}).
\end{remark}

\smallskip

\begin{proposition}[Variation of $F$]\label{prop:DeltaF33}
  Consider the interactions of any two wave fronts of the same family,
  $1$ or $3$, and assume (\ref{eq:Ao>0}), (\ref{rogna}). If
\begin{eqnarray}\label{eq:sogliazza}
1 < \xi < \frac{1}{d} \qquad {and}\qquad K < \frac{\xi-1}{A_o}
\end{eqnarray}
then $\Delta L_{\xi} < 0$ and $ \Delta F < 0$\,.
\end{proposition}

\smallskip

\begin{proof}
  Let two waves $\alpha_i$, $\beta_i$ interact, $i=1,3$, giving rise
  to waves $\eps_1$, $\eps_3$. We consider $i=3$, the other case being
  analogous.
Using (\ref{eq:333}), we get
\begin{eqnarray}\label{Delta-Q-interaz-13}
\Delta Q = \left(|\eps_3| - |\alpha_3| -
    |\beta_3|\right) L_{cd}^+ + |\eps_1|L_{cd}^- & \leq &  |\eps_1|L_{cd}^-  \,\leq\, |\eps_1| A_o \,.
\end{eqnarray}
Now we claim that
\begin{eqnarray}\label{Delta_L_xi_13}
\Delta L_{\xi} + |\eps_1|(\xi - 1) & \le & 0\,.
\end{eqnarray}
From this estimate it follows $\Delta F = \Delta L_{\xi} + K
\Delta Q \le |\eps_1|\left(1-\xi+K A_o\right)< 0$ because of
(\ref{eq:sogliazza}). To prove our claim we consider the possible
cases; we make use of (\ref{1-0}).

\smallskip
\fbox{$SS\to RS$}\quad
Since $\Delta L=0$ then
\begin{eqnarray}\label{Delta_xi_L_xi_SSRS}
\Delta L_{\xi} + (\xi -1) |\eps_1| = \xi (|\eps_1| +
|\eps_3|-|\alpha_3|- |\beta_3|)\le 0\,.
\end{eqnarray}

\smallskip

\fbox{$SR,\, RS\to SR$}\quad Assume $\alpha_3<0<\beta_3$; then
(\ref{Delta_L_xi_13}) reads $(2\xi-1) |\eps_1| + |\eps_3|-
|\beta_3|-\xi|\alpha_3|\le 0$. For later use we prove the stronger
inequality
\begin{eqnarray}\label{Delta_xi_L_xi_SRSR}
\xi^2 |\eps_1| + |\eps_3|- |\beta_3|-\xi|\alpha_3|\le 0\,.
\end{eqnarray}
Indeed, from Lemma~\ref{lem:shock-riflesso}{\it (ii)} we have
$|\eps_3|< |\beta_3|$, while $\xi |\eps_1|\le |\alpha_3|$ from
(\ref{eq:chi_def}), (\ref{eq:sogliazza})${}_1$.

\smallskip \fbox{$SR,\, RS\to SS$}\quad Assume
$\alpha_3<0<\beta_3$; then (\ref{Delta_L_xi_13}) is $(2\xi-1)
|\eps_1| + \xi(|\eps_3|-|\alpha_3|) - |\beta_3|\le 0$. We prove
also in this case the stronger inequality
\begin{eqnarray}\label{Delta_xi_L_xi_SRSS}
\xi^2 |\eps_1| + \xi(|\eps_3|-|\alpha_3|) - |\beta_3|\le 0.
\end{eqnarray}
Indeed, by (\ref{1-0s}) and again because of (\ref{eq:chi_def}),
(\ref{eq:sogliazza})${}_1$, one has
\begin{eqnarray*}
\xi^2 |\eps_1| + \xi(|\eps_3|-|\alpha_3|) - |\beta_3|&=& \xi^2 |\eps_1| + \xi(|\eps_1|-|\beta_3|) - |\beta_3|\\
&=& (\xi+1) (\xi|\eps_1| - |\beta_3|) \le 0\,.
\end{eqnarray*}

This proves the claim and concludes the proof.
\end{proof}

\smallskip

\begin{remark}\rm
From the above proof we see that $\Delta L_\xi\le 0$ for $\xi=1$.
This was a key point in \cite{Nishida68}, where however a
different choice of strengths was done. In \cite{AmadoriGuerra01}
the inequality $\Delta L_\xi\le 0$ was proved to hold also for
$1<\xi\leq \xi_o$, for some $\xi_o>1$; the condition
(\ref{eq:sogliazza})${}_1$ gives an estimate of such a threshold.

More precisely, in the first two cases of Proposition
\ref{prop:DeltaF33} we have $\Delta L_\xi\le 0$ for every
$\xi\ge1$. The third case is analyzed in detail in Lemma
\ref{lem:second}; we prove there that $\Delta L_\xi\le 0$ for any
$\xi > 1$ if $c(m) \le 1/2$, while we need $1<\xi\le
\frac{1}{2c(m)-1}$ if $c(m)>1/2$.
\end{remark}


\subsection{Non-physical waves}
In this subsection we compute the strength of a non-physical wave
generated by an interaction and prove that it does not change in
subsequent interactions. We introduce the following notation: given
$U_\ell=(v_\ell,u_\ell,\lambda_\ell)$ and $\lambda_r$ we define by
\begin{eqnarray*}
U_{\ell r}^* = \Phi_2(\lambda_r,U_\ell) = \left(A_{r\ell}
v_\ell,u_\ell, \lambda_r \right)
\end{eqnarray*}
the state on the right of a $2$-wave with left state
$U_\ell=(v_\ell,u_\ell,\lambda_\ell)$ and $\lambda=\lambda_r$ on the
right, where $A_{r\ell} = a^2(\lambda_r) / a^2(\lambda_\ell)$. See
(\ref{eq:lax2}) and \cite{AmCo06-Proceed-Lyon}.


\begin{figure}[htbp]
\begin{picture}(100,100)(-80,-15)
\setlength{\unitlength}{1pt}

\put(20,0){
\put(0,0){\line(0,1){40}} \put(-2,-5){\makebox(0,0){$\delta_2$}}
\put(0,40){\line(1,-1){30}} \put(38,5){\makebox(0,0){$\delta_1$}}
\put(0,40){\line(0,1){30}} \put(-2,78){\makebox(0,0){$\delta_2$}}
\put(0,40){\line(5,1){50}} \put(60,50){\makebox(0,0){$n\!p$}}
\put(0,40){\line(-1,1){30}}
\put(-30,77){\makebox(0,0){$\delta_1$}}
\put(-25,40){\makebox(0,0){$U_\ell$}}
\put(14,11){\makebox(0,0){$U_{\ell r}^*$}}
\put(40,30){\makebox(0,0){$U_r$}}
\put(25,60){\makebox(0,0){$U_{qr}^*$}}
\put(-8,60){\makebox(0,0){$U_q$}}
\put(0,-20){\makebox(0,0){$(a)$}}
}

\put(190,0){
\put(0,0){\line(0,1){40}}
\put(-2,-5){\makebox(0,0){$\delta_2$}}
\put(0,40){\line(-1,-1){30}}
\put(38,75){\makebox(0,0){$\delta_3$}}
\put(0,40){\line(0,1){30}}
\put(-2,78){\makebox(0,0){$\delta_2$}}
\put(0,40){\line(5,1){50}}
\put(60,50){\makebox(0,0){$n\!p$}}
\put(0,40){\line(1,1){30}}
\put(-30,3){\makebox(0,0){$\delta_3$}}
\put(-25,40){\makebox(0,0){$U_\ell$}}
\put(40,20){\makebox(0,0){$U_{nr}^*=U_r$}}
\put(40,60){\makebox(0,0){$U_q$}}
\put(12,65){\makebox(0,0){$U_{lr}^*$}}
\put(-10,10){\makebox(0,0){$U_n$}}
\put(0,-20){\makebox(0,0){$(b)$}}
}

\end{picture}

\caption{\label{fig:interNP}{Simplified Riemann solver
}}
\end{figure}
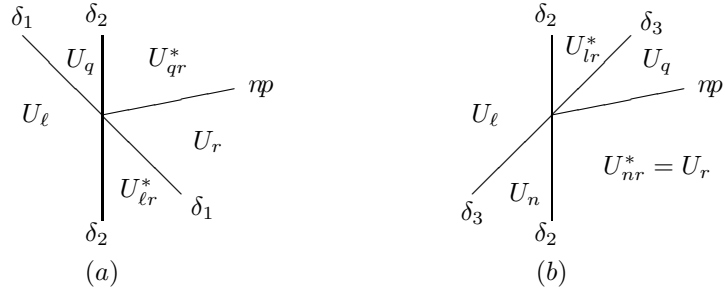


\smallskip

\begin{proposition}[Non-physical waves]\label{lem:np}
Consider $U_\ell=(v_\ell,u_\ell,\lambda_\ell)$. Let
$U_r=(v_r,u_r,\lambda_r)$ be connected to $U^*_{\ell r}$ by a
1-wave of size $\delta_1$ and $U_q=(v_q,u_q,\lambda_\ell)$ be
connected to $U_\ell$ by a 1-wave of size $\delta_1$, see Figure
\ref{fig:interNP}(a). Assume (\ref{rogna}).

Then $U^*_{qr}$ and $U_r$ differ only in the $u$ component; if
$\delta_2$ denotes the size of the 2-wave, there exists a constant
$C_o=C_o(m)$
such that
\begin{eqnarray}\label{size-onda-nonfisica}
\| U^*_{qr} - U_r\| = |u_q - u_r| \leq C_o |\delta_2 \delta_1|\,.
\end{eqnarray}
A similar result holds for the interaction of a 3-wave, see Figure
\ref{fig:interNP}(b), again under (\ref{rogna}).

Moreover the size of a non-physical wave does not change in
subsequent interactions. For any $K>0$ and $\knp < K/C_o$ at any
interaction involving a non-physical wave we have
\begin{equation}\label{eq:FdecrNP}
\Delta F  \le 0\,,
\end{equation}
with $\Delta F  < 0$  when a non-physical wave is generated.
\end{proposition}

\smallskip

\begin{proof}
Recalling \cite[Lemma 2]{AmCo06-Proceed-Lyon}, only the $u$
component will be different after commutation of the $1$- and the
$2$-wave. We find that
$$
u_q - u_\ell = 2a_\ell h(\delta_1),\qquad u_r - u_\ell = 2a_r
h(\delta_1)
$$
hence
$$
|u_q - u_r| = 2|a_\ell - a_r| \cdot |h(\delta_1)|\leq |\delta_2
\delta_1| \cdot 2 a(1) \max_{0<\eta\leq m} \frac{\sinh
\eta}\eta\,.
$$
Then (\ref{size-onda-nonfisica}) follows with $C_o(m) \doteq 2
a(1)\cdot \frac{\sinh m}{m}$.

\smallskip

Next, assume that a non-physical wave interacts with a $2$-wave.
Since the values of $u$ do not change across a $2$-wave, the left
and right values of $u$ of the non-physical wave do not change
across the interaction; hence the size does not change.

Assume then that a non-physical wave interacts with a $1$- or
$3$-wave of size $\delta$. Since $\lambda$ is constant, we refer
only to the components $v$, $u$. Let $(v_\ell,u_\ell)$ and
$(v_\ell,u_q)$ be the side states of the non-physical wave before
the interaction and $(v_\ell,u_q)$, $(v_r,u_r)$ be the side states
of the physical wave. After the interaction, let $(\tilde
v_\ell,\tilde u_\ell)$ be the intermediate state. One has
$$
u_r - u_q = 2 a(\lambda)h(\delta) = \tilde u_\ell - u_\ell\,,
$$
hence $|u_q - u_\ell| = |u_r-\tilde u_\ell|$.

At last, we consider the functional $F$. The potential $Q$ is
unaltered when non-physical waves interact with other waves. The
only cases in which $L_{np}$ changes are when a non-physical wave
arises. Assume that a $1$- or a $3$-wave of size $\delta$
interacts with a $2$-wave of size $\delta_2$, producing a wave of
same size and a non-physical wave. Then $\Delta Q = -
|\delta\delta_2|$ and $\Delta L_\xi = \knp \Delta L_{np} \leq \knp
C_o|\delta\delta_2|$; hence $\Delta F = \Delta L_\xi + K \Delta Q
\leq |\delta\delta_2| (\knp C_o - K)$, then (\ref{eq:FdecrNP}).
\end{proof}


\subsection{Decreasing of the functional $F$ and control of the variations}
We collect first the previous results into a single proposition.

\smallskip

\begin{proposition}[Local decreasing]\label{prop:last}
  Consider the interaction of any two waves either of families $1$, $2$, $3$
  or non-physical. Assume (\ref{rogna}) for some $m > 0$;
  let $C_o=C_o(m)$ as in Proposition \ref{lem:np}. Finally let $A_o$
  satisfy
\begin{eqnarray}\label{bound_su_A_o}
0 < A_o < 2\frac{1-d}{3-d}\,.
\end{eqnarray}
If $\xi$, $K$, $\knp$ satisfy
\begin{eqnarray}\label{cond_su_chi_xi}
\frac{2-A_o}{2-3A_o} < \xi <\frac{1}{d}\,,\qquad \qquad
\frac{2\xi}{2-A_o}<K<\frac{\xi-1}{A_o}\,,\qquad \qquad
\knp<\frac{K}{C_o}
\end{eqnarray}
then
\begin{equation}\label{eq:Fdecr}
\Delta F \le 0\,.
\end{equation}
\end{proposition}

\begin{proof}
The condition on $K$ comes from (\ref{xi_e_K}) and
(\ref{eq:sogliazza})${}_2$. The interval where $K$ lies is not
empty if $A_o<\frac23$ and $\xi > \frac{2-A_o}{2-3A_o}$; together
with (\ref{eq:sogliazza})${}_1$ this gives the assumption required
on $\xi$. In turn, it is possible to choose $\xi$ in such interval
if (\ref{bound_su_A_o}) holds. Remark that $2\frac{1-d}{3-d} \le
\frac23$, so the previous condition on $A_o$ holds. Therefore the
assumptions of Propositions~\ref{Delta-F-2wave},
\ref{prop:DeltaF33} and \ref{lem:np} hold and then
(\ref{eq:Fdecr}) follows.
\end{proof}

\smallskip

\begin{proposition}[Global decreasing] \label{prop:very_last}
Let $m > 0$; assume (\ref{bound_su_A_o}), (\ref{cond_su_chi_xi}),
\begin{eqnarray}
L(0+) < \frac{m}{2\xi-1}\label{eq:super_cond}
\end{eqnarray}
and that the approximate solution $U$ is defined in $[0,T]$. Then
$L(t) < m$ for any $t\in[0,T]$; as a consequence, condition
(\ref{rogna}) holds for any $1$ or $3$ wave in $U$. Finally,
$\Delta F(t)\le0$ for all times $t\in[0,T]$.
\end{proposition}

\smallskip

\begin{proof}
  Because of Proposition \ref{prop:last}, for any time $t$ not of
  interaction we have
\[
L(t) \le L_\xi(t) \le F(t) \le F(0+) \le \xi L(0+)+KQ(0+)\,.
\]
Moreover $Q(0+)\le L(0+)L_{cd} \le L(0+)A_o$; thanks to
(\ref{eq:super_cond}) and (\ref{cond_su_chi_xi}) we have
\begin{equation}
L(t) ~ \le ~ L(0+) \cdot \left(\xi + K A_o\right) ~ \le ~ L(0+)
\cdot \left(2\xi -1\right)~ < ~ m\,. \label{eq:opt}
\end{equation}
Therefore the size $|\epsilon|$ of any wave of families $1$, $3$
at time $t>0$ satisfies $|\eps|\leq m$, so (\ref{rogna}) holds.
\end{proof}

\begin{remark}\label{rem:bestchoice}\rm From (\ref{eq:super_cond}) we see that,
in order to have $L(t) < m$, the smaller is $\xi$ the larger can
be chosen $L(0+)$.
\end{remark}


\section{The convergence and the consistence of the
  algorithm}\label{sec:Cauchy} In this section we prove Theorem
\ref{thm:main}. We show first that for fixed $\nu$ the algorithm
introduced in Section \ref{sec:app_sol} gives an approximate solution
defined for every $t>0$; more precisely we prove that at every time
the number of interactions is bounded. Then we prove that the total
amount of non-physical waves in each approximate solution is very
small. The convergence of a suitable subsequence is assured by Helly's
theorem; then consistence follows.


\subsection{Control of the number of interactions}\label{subs:number_int}
We prove first that the size of the rarefactions in the scheme is
small.
\begin{lemma}\label{lem:small_rar}
Consider a rarefaction of size $\eps$; then
\begin{eqnarray}\label{size-raref}
|\eps|< \eta{\rm e}^{\frac{A_o}2}\,.
\end{eqnarray}
\end{lemma}

\begin{proof} We analyze all possible situations. When the rarefaction
  is generated, one has $0<\eps<\eta$. When it interacts with a 1- or
  3-wave, the size does not increase, see Remark
  \ref{rem:rar_decreasing}. By Proposition~\ref{lem:np} the size does not
  change when interactions with non-physical waves occur.

  The last case to be considered is when a rarefaction interacts with
  a 2-wave. In this case the size may increase; however, a rarefaction
  can meet a fixed 2-wave only once. Consider the case of a
  1-rarefaction of size $\delta_1$, as in
  Proposition~\ref{Delta-F-2wave}, the other being analogous. If
  $\delta_2<0$ then the size decreases, see (\ref{eq:12-}). If
  $\delta_2>0$ by (\ref{eq:12+}) we have
$$
|\eps_1| = |\delta_1| +  |\eps_3|\leq |\delta_1|\left( 1 + \frac 1 2
  |\delta_2|\right) < |\delta_1| {\rm e}^{\frac{|\delta_2|}2}\,.
$$
Summarizing the three cases above, we get $|\eps|< \eta{\rm
  e}^{L_{cd}^+/2}$ (or $|\eps|< \eta{\rm e}^{L_{cd}^-/2}$) for a
$1$-rarefaction (resp. $3$-rarefaction), where $L_{cd}^\pm$ is the
sum of the $2$-waves at the right or left of the rarefaction. Then
(\ref{size-raref}) follows.
\end{proof}

Next, we prove that the number of interactions remains bounded in
finite time, so that the approximate solution is well defined for all
$t>0$. We give first a lemma.

\smallskip

\begin{lemma}\label{lem:wfai-ii} Consider the wave-front tracking algorithm described in
  Section \ref{sec:app_sol}, under the assumptions of
  Proposition~\ref{prop:very_last}. Then

\begin{romannum}
\item the number of interactions involving a $2$-wave and
solved by the accurate Riemann solver is finite;

\item the number of interactions where a new rarefaction of
size $\eps\geq \eta$ arises is finite.
\end{romannum}
\end{lemma}

\smallskip

\begin{proof}
  Consider first \emph{(i)}
  and refer to Proposition~\ref{Delta-F-2wave}. Then, using
  (\ref{diseq:F-xi-K}), we have $\Delta F\le
  \rho\left(\xi+K(A_o-2)/2\right)<0$, hence $F$ decreases by a
  uniform positive quantity; since it is non increasing, this can
  happen only a finite number of times.

\smallskip

Consider then \emph{(ii)}. After \emph{(i)}, it remains only to
consider the case of two shocks of the same family interacting.

Under the notation of the corresponding case in the proof of
Proposition~\ref{prop:DeltaF33} we have $\eps=\eps_1\geq \eta$ and
\begin{eqnarray*}
\Delta F \le |\eps_1|\left(1-\xi+K A_o\right) &\le& \eta
\left(1-\xi+K A_o\right)<0
\end{eqnarray*}
because of (\ref{cond_su_chi_xi}). Arguing as in \emph{(i)}, this
can happen only a finite number of times.
\end{proof}

\smallskip

About \emph{(ii)} in Lemma \ref{lem:wfai-ii}, recall that if $\eps \ge \eta$ then
the new rarefaction must be split into more than one wave.
Therefore Lemma \ref{lem:wfai-ii} can be rephrased by saying that,
except for finite interactions, the number of waves emitted in an
interaction is at most three, and this case occurs precisely when
a non-physical wave is generated; moreover, in every interaction
at most one wave per family is emitted.

In a schematic way, apart from a finite number of interactions, in our
algorithm the following holds (we will consider the set of
non-physical waves as a $4^{th}$ family of waves):

\smallskip

\begin{itemize}

\item[(a)] the interaction of an $i$-wave, $i=1,3$, with a $2$-wave is
solved by a single $i$-wave, a $2$-wave and a $4$-wave;

\item[(b)] in the interaction of just $1$- and/or $3$-waves, there is
at most one outgoing wave of each family $1$ and $3$;

\item[(c)] the interaction of a $i$-wave, $i=1,2,3$, with a $4$-wave is
solved by an $i$-wave and a $4$-wave.

\end{itemize}

\smallskip

\noindent The following proposition is inspired by \cite[Lemma
2.5]{AmadoriGuerra01}.

\smallskip

\begin{proposition}\label{prop:finite_interazioni}
  Consider the wave-front algorithm described in Section
  \ref{sec:app_sol} and assume in the strip $[0,T)\times \reali$ the
  following:
\begin{quote}
  for some $a_1<a_2<0<b_1<b_2$ the waves of the first (third) family
  have speeds in the interval $[a_1,a_2]$ (resp. $[b_1,b_2]$).
\end{quote}
Then the number of interactions in the region $[0,T)\times \reali$ is
finite.
\end{proposition}

\smallskip

\begin{proof}
  Assume by contradiction that in the region $[0,T)\times \reali$
  there exists an infinite number of interactions. Unless of taking a
  smaller $T$ we can assume that the number of interactions is finite
  in every strip $[0,t]\times\reali$, $0<t<T$, and that $T$ is an
  accumulation point for the times of interaction. Then there exists a
  sequence $(t_j,x_j)$, $j=1,2,\ldots$ of interaction points such that
\[
0< t_j < t_{j+1} < T\  \hbox{ for all } j  \ \hbox{ and }
(t_j,x_j)\to (T,\bar x)
\]
for some $\bar x$. Denote ${\cal J} = \left\{(t_j,x_j)\colon
  j=1,2,\ldots\right\}$ the set of all interaction points; the set
${\cal J}$ is bounded in the strip $[0,T)\times \reali$ because of the
finite propagation speed.

The situation described in items (a)--(c) above holds except a
finite number of interactions; let $\tau < T$ be the maximum time
of these ``exceptional'' interactions. It is not restrictive to
assume that $\tau<t_j<T$ for all $j=1,2\ldots$.

Starting from a point of ${\cal J}$, we ``trace back'' all the
segments up to $t=0$; we repeat the procedure for all the points of
${\cal J}$ and call ${\cal F}$ the set of the ``traced segments''
obtained in this way. In other words, a segment belongs to the set
${\cal F}$ iff it can be joined forward in time to some point of
${\cal J}$ by a continuous path along the wave fronts.  The set ${\cal
  F}$ is not empty: for instance, two segments interacting at the
point $(t_j,x_j)$ belong to ${\cal F}$, for any $j =1,2,\ldots$.
Observe the following dichotomy property of $\mathcal{F}$ which is
used just below:
\begin{quote}
  two interacting waves either belong both to ${\cal F}$ or none of
  them does; moreover if at least one of the outgoing waves belong to
  ${\cal F}$, then both the incoming waves must belong to ${\cal F}$.
\end{quote}
We partition now all the interaction points of the algorithm that
occur for times $t>\tau$ into the following sets:

\smallskip

\begin{itemize}

\item ${\cal I}_0$: the interaction points where no ingoing wave
belongs to ${\cal F}$;

\item ${\cal I}_1$: the interaction points where both incoming
waves belong to $\mathcal{F}$ and at most one outgoing segment
belongs to ${\cal F}$;

\item ${\cal I}_2$: the interaction points where exactly two
outgoing segments belong to ${\cal F}$;

\item ${\cal I}_3$: the interaction points where three outgoing
waves all belong to ${\cal F}$.

\end{itemize}

\smallskip

\noindent Because of the dichotomy property quoted above no outgoing wave in
case ${\cal I}_0$ can belong to $\mathcal{F}$. On the contrary,
both incoming waves in $\mathcal{I}_1$, $\mathcal{I}_2$,
$\mathcal{I}_3$ must belong to $\mathcal{F}$.

Recall that we are considering times $t > \tau$. Therefore the
maximum number of emitted waves in an interaction is three, and
this happens only in the situation considered in $\mathcal{I}_3$,
that is for interactions as in \emph{(a)} above. The case of more
than one emitted wave per family cannot occur, and so the outgoing
waves in ${\cal I}_2$ belong to different families. By definition
we have ${\cal J}\cap{\cal I}_0=\emptyset$ and so
$\mathcal{J}\subset \mathcal{I}_1 \cup \mathcal{I}_2 \cup
\mathcal{I}_3$.

Let ${\cal V}(t)$ be the total number of wave-fronts of the
families $1$, $2$ and $3$ that belong to ${\cal F}$, at time $t$.
The functional ${\cal V}(t)$ is non-increasing, and it decreases
at least by $1$ across ${\cal I}_1$. Then ${\cal I}_1$ is finite.
As a consequence, all the interaction points of $\mathcal{J}$
belong to ${\cal I}_2\cup{\cal I}_3$, except at most a finite
number. Let $\tau_1\in\left[\tau, T\right)$ be a time such that
all points in ${\cal I}_1$ lies in $t < \tau_1$.

\smallskip

Let $P=\{x_1,\ldots,x_{N_1}\}$ the set of points of the $x$-axis where
a 2-wave is located, $N_1\leq N_o$. We consider two cases.

\smallskip

\paragraph{\fbox{$\bar x\not\in P$}} In this case we can choose the time $\tau_1<T$
such that after that time no segment belonging to ${\cal F}$ crosses a
2-wave. Then all the points in ${\cal J}$ with $t_j>\tau_1$ belong to
${\cal I}_2$. The same argument of \cite[Lemma 2.5]{AmadoriGuerra01}
can now be used, reaching a contradiction.

\smallskip

\paragraph{\fbox{$\bar x\in P$}} Consider $(t^*,x^*)\in {\cal J}$ with $t^* > \tau_1$.
As in \cite[Lemma 2.5]{AmadoriGuerra01} we define for $t\in\left[
t^*, T \right)$ two continuous paths $\gamma_\ell(t)$,
$\gamma_r(t)$ starting at $(t^*,x^*)$ in the following way.

At $(t^*,x^*)$ there are either two or three outgoing segments
belonging to ${\cal F}$; for times $t>t^*$ and sufficiently close to
$t^*$ we define $\gamma_\ell(t)$ to be the segment on the left and
$\gamma_r(t)$ the one on the right. When $\gamma_\ell(t)$
($\gamma_r(t)$) arrives at an interaction point, it is prolonged by a
segment of ${\cal F}$; since they are at least two, it follows the one
on the left (resp., on the right). Then the path $\gamma_\ell(t)$ is
made by segment of the families $1$, $2$ or $3$, while $\gamma_r(t)$
is made by segments of families $3$ or $4$; in fact a look to the
proof of Proposition~\ref{lem:np} shows that the interaction of a
$1$-wave with a $2$-wave always produces a $4$-wave.

We claim that the speeds of the paths $\gamma_\ell(t)$, $\gamma_r(t)$
are strictly separated. This is clear for a finite number of nodes
among segments. If it happens that $\gamma_r(t)$ follows once a
$4$-segment, then from that point on it will follow only $4$-segments,
so $\dot\gamma_\ell(t)\le b_2 < \hat{s} = \dot\gamma_r(t)$. If on the
other hand $\gamma_r(t)$ never follows a $4$-segment, then
$\dot\gamma_\ell(t)\le 0 < b_1 \le \dot\gamma_r(t)$.

This proves the claim; the same argument of \cite{AmadoriGuerra01}
then leads to a contradiction.
\end{proof}


\subsection{Control of the total size of non-physical
  fronts}\label{subs:control} Assume as above that the assumptions of
Proposition \ref{prop:very_last} hold. We assign inductively to each
wave $\alpha$ a \emph{generation order} $k_\alpha$ as in \cite[p.
140]{Bressanbook}. This is done according to the following procedure.
First, at time $t=0$ each wave has order~$1$. Second, assume that two
waves $\alpha$ and $\beta$ interact at time $t$; if $\alpha$ and
$\beta$ belong to different families, the outgoing waves of those
families keep the order of the incoming waves, the other waves assume
order $\max\{k_\alpha, k_\beta\}+1$; if $\alpha$ and $\beta$ belong to
the same family, the outgoing wave of that family takes the order
$\min\{k_\alpha, k_\beta\}$, the other waves are assigned order
$\max\{k_\alpha, k_\beta\}+1$.

When specialized to the current setting this has the following
consequences:
\begin{itemize}
\item every $2$-wave has order $1$; when a $i$-wave, $i=1,3$, of order
  $k$ interacts with a $2$-wave, the outgoing $i$-wave has order $k$,
  the other outgoing wave (of the family $j$, $j=1,3$, $j\ne i$, or a
  non-physical wave) has order $k+1$;
\item in the interaction of a $1$- with a $3$-wave the waves cross
  without changing order; in the interaction of two waves $\alpha$,
  $\beta$ of the same family $i=1,3$, the outgoing wave of the family
  $i$ takes order $\min\{k_\alpha, k_\beta\}$, the wave of the family
  $j=1,3$, $j\ne i$, has order $\max\{k_\alpha, k_\beta\}+1$;
\item when a non-physical wave interacts with any other wave, both
  waves cross without changing order; in particular a non-physical
  wave keeps the order it has been assigned when generated.
\end{itemize}

\noindent For $t\ge0$ not an interaction time and any
$k=1,2,\ldots$ define, see (\ref{L-xi}), (\ref{Q}),
\begin{eqnarray*}
V_k(t)  &=&  \sum_{\gamma>0\atop k_\gamma = k}|\gamma| + \xi
\sum_{\gamma<0\atop k_\gamma = k}|\gamma| +
\knp \sum_{\gamma \in {\cal N\!P}\atop k_\gamma = k}|\gamma|
\\
Q_k(t) &=& \sum_{\gamma_3 \mbox{\sz \ at the left of }
\delta_2\atop k_{\gamma_3} = k} |\gamma_3| |\delta_2| + \sum_{\gamma_1
\mbox{\sz\ at the right of } \delta_2 \atop k_{\gamma_1} = k}
|\delta_2| |\gamma_1|
\\
F_k(t) &=& V_k(t) + K Q_{k}(t)
\end{eqnarray*}
and
\begin{eqnarray*}
\tilde V_k(t) =
 \sum_{\ell\geq k} V_\ell(t)\,,\qquad \tilde Q_k(t) =
 \sum_{\ell\geq k} Q_\ell(t)\,,\qquad
\tilde F_k(t) = \tilde V_k(t) + K\tilde Q_{k}(t)\,.
\end{eqnarray*}
We remark that
\begin{eqnarray*}
\tilde F_1(0+)=L_\xi(0+)+K Q(0+),\qquad \tilde F_k(0+)=0\qquad \mbox{for } k\ge 2.
\end{eqnarray*}
Observe that if a non-physical front interacts with another wave, the
functionals above do not change; the same holds for interactions
between $3$- and $1$-waves. Then we focus on interactions of waves of
the same $i$ family, $i=1,3$ (as usual denote $j=1,3$, $j\ne i$) and
on interactions between $1$- or $3$-waves with a $2$-wave.

For $h\in\naturali$, denote by $I_h$ the set of times $t$ when an
interaction occurs between two waves $\alpha$ and $\beta$ of families
$1$ or $3$ with $\max\{k_\alpha, k_\beta\}=h$; denote by $J_h$ the set
of interaction times $t$ of a $1$- or $3$-wave of order $h$ with a
$2$-wave.  Finally, denote ${\cal T}_h = I_h\cup J_h$ and
$I=\bigcup_{h\ge 1} I_h$, $J=\bigcup_{h\ge 1} J_h$, ${\cal T}=I\cup
J$.

In order to control the total size of non-physical fronts we must
strengthen the assumptions (\ref{bound_su_A_o}) and
(\ref{cond_su_chi_xi}) required in Proposition
\ref{prop:very_last}. First, for any fixed $m>0$, instead of
(\ref{bound_su_A_o}) we require the stronger condition
\begin{eqnarray}
0 ~ < ~ A_o & < & \frac{1-\sqrt{d}}{2-\sqrt{d}}\,.
\label{bound_su_A_o-NEW}
\end{eqnarray}
Then denote
\begin{equation}
\lambda ~ \dot = ~ \frac{1+KA_o}{\xi}\,, \qquad\lambda_2 ~ \dot =
~ \frac{\xi+KA_o}{K(2-A_o) - \xi}\,, \qquad \mu ~ \dot = ~
\max\left\{\lambda \,,\lambda_2, \frac{\knp C_o}{K}\right\} \,.
\label{llk}
\end{equation}
We need $0<\mu<1$. From (\ref{cond_su_chi_xi})${}_3$ we have $\knp
C_o/K<1$; moreover we have $0<\lambda<1$ and $\lambda_2>0$ because
of (\ref{cond_su_chi_xi})${}_2$. At last $\lambda_2<1$ holds iff
$\frac{\xi}{1-A_o}<K$; this condition is stronger than the left
inequality in (\ref{cond_su_chi_xi})${}_2$. Hence, instead of
(\ref{cond_su_chi_xi}) we assume
\begin{eqnarray}\label{cond_su_chi_xi-NEW}
\frac{1-A_o}{1-2A_o} < \xi < \frac{1}{\sqrt{d}} \,,\qquad \qquad
\frac{\xi}{1-A_o}<K<\frac{\xi-1}{A_o}\,,\qquad \qquad \knp <
\frac{K}{C_o}\,.
\end{eqnarray}
Remark the new upper bound required on  $\xi$. As in the proof of
Proposition \ref{prop:last}, the interval where $K$ varies is not
empty if $A_o < \frac12$ and $\xi >\frac{1-A_o}{1-2A_o}$. In turn,
we can find $\xi$ satisfying (\ref{cond_su_chi_xi-NEW})${}_1$ if
(\ref{bound_su_A_o-NEW}) holds; remark that $
\frac{1-\sqrt{d}}{2-\sqrt{d}}\le \frac12$. The last condition in
(\ref{cond_su_chi_xi-NEW}) coincides with that in
(\ref{cond_su_chi_xi}).

\smallskip

\begin{remark}\label{rem:bestbest}\rm
Proposition \ref{prop:very_last} still holds under the stronger
assumptions (\ref{bound_su_A_o-NEW}), (\ref{cond_su_chi_xi-NEW})
under the same condition (\ref{eq:super_cond}), because the
inequality on the right hand side in
(\ref{cond_su_chi_xi-NEW})${}_2$ has not changed, see
(\ref{eq:opt}).
\end{remark}

\smallskip


\begin{proposition}\label{lem:01}
Fix $m>0$ and assume (\ref{bound_su_A_o-NEW}),
(\ref{cond_su_chi_xi-NEW}). We have:
\begin{remunerate}
\item $\tau\in {\cal T}_h$, $h\le k-2$: then $\Delta\tilde F_k=\Delta
  F_k=0$\,.
\item $\tau\in {\cal T}_{k-1}$: then $\Delta F_{k-1}<0$,
$\Delta\tilde F_k=\Delta F_k>0$
and
\begin{eqnarray}\label{h=k-1b-XXX}
\piu{\Delta\tilde F_k} \le \mu \Bigl(\meno{\Delta F_{k-1}}   -
\sum_{\ell=1}^{k-2} \piu{\Delta F_\ell} \Bigr)\,.
\end{eqnarray}
\item $\tau\in {\cal T}_h$, $h\ge k$: if $h=k$ then $\Delta F_k < 0$;
  in any case $\Delta \tilde F_k < 0$ and
\begin{eqnarray}\label{Delta_tilde_V_k-meno}
\sum_{\ell=1}^{k-1} \piu{\Delta F_{\ell}}
& < & \meno{\Delta\tilde F_k}\,.
\end{eqnarray}
\end{remunerate}
\end{proposition}


\smallskip

\begin{proof}
  As we pointed out above, for $\tau \in I$ only interactions of waves
  of the same family are taken into account. Remark now that for
  interactions between $3$-waves we have $\Delta Q_\ell = L_{cd}^+ \Delta
  V_\ell$ while $\Delta Q_\ell = L_{cd}^- \Delta V_\ell$ holds for
  1-waves. So if $\tau \in I$ then $\Delta V_\ell>0$ ($\Delta
  V_\ell<0$, $\Delta V_\ell=0$) implies $\Delta Q_\ell\ge 0$ (resp.
  $\Delta Q_\ell\le 0$, $\Delta Q_\ell= 0$).

  Remark also that by Proposition~\ref{prop:very_last}
\begin{eqnarray}\label{Delta_W_piu_V}
\sum_{\ell=1}^{k-1}\Delta F_{\ell} ~+~ \Delta
\tilde F_k <0 \,.
\end{eqnarray}
%

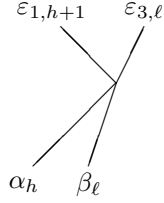
\begin{figure}[htbp]
\begin{picture}(100,60)(-80,-10)
\setlength{\unitlength}{0.7pt}

\put(180,0){ \put(0,40){\line(-1,-1){45}}
\put(-50,-15){\makebox(0,0){$\alpha_{h}$}}
\put(0,40){\line(-1,-3){15}}
\put(-15,-15){\makebox(0,0){$\beta_{\ell}$}}
\put(0,40){\line(1,2){15}}
\put(15,78){\makebox(0,0){$\eps_{3,\ell}$}}
\put(0,40){\line(-1,1){30}}
\put(-35,78){\makebox(0,0){$\eps_{1,h+1}$}}
}

\end{picture}

\caption{\label{fig:inter_k-1}{Interactions of $3$-waves;
$h\ge\ell$ denote generation orders.}}
\end{figure}

\begin{remunerate}

\item If $h\le k-2$, no waves with order $\ge k$ are involved
  and then $\Delta\tilde F_k=\Delta F_k=0$\,.

\item Let $h=k-1$. First, consider $\tau\in I_{k-1}$; then $\Delta
  \tilde V_k=\Delta V_k>0$. We prove that
\begin{eqnarray}\label{0k-1}
\piu{\Delta \tilde V_k} & \le & \frac{1}{\xi} \Bigl(\meno{\Delta
V_{k-1}} - \sum_{\ell = 1}^{k-2}\piu{\Delta V_\ell} \Bigr)\,.
\end{eqnarray}
Indeed, from (\ref{Delta_xi_L_xi_SSRS})--(\ref{Delta_xi_L_xi_SRSS}),
we deduce
\begin{eqnarray}\label{eq:ok}
\xi\piu{\Delta \tilde V_k} + \Delta V_{k-1} +
\sum_{\ell=1}^{k-2}\Delta V_{\ell} \le 0\,.
\end{eqnarray}
If $\min\{k_\alpha,k_\beta\}=k-1$, then $\Delta V_\ell =0$ for
$\ell=1,\ldots,k-2$ and (\ref{eq:ok}) becomes $\xi[\Delta \tilde
V_k]_+ + \Delta V_{k-1}<0$; this implies $\Delta V_{k-1}<0$ and hence
$[\Delta \tilde V_k]_+ < (1/\xi)[\Delta V_{k-1}]_-$, that is
(\ref{0k-1}).

If $\min\{k_\alpha,k_\beta\}=\ell\le k-2$, then $\Delta V_{k-1}<0$
since no waves of order $k-1$ are present after the interaction;
therefore the estimate (\ref{eq:ok}) becomes
$$
\xi\piu{\Delta \tilde V_k} - \meno{\Delta V_{k-1}} +\Delta
V_{\ell} \le 0\,.
$$
If $\Delta V_{\ell}\ge 0$ we get (\ref{0k-1}). If $\Delta
V_{\ell}<0$ we check directly that $\xi[\Delta \tilde V_k]_+ -
[\Delta V_{k-1}]_-\le \xi^2 |\eps_1| - |\alpha_3| <0$ because of
(\ref{cond_su_chi_xi-NEW})${}_1$. This proves (\ref{0k-1}); from
this proof we see moreover that $\Delta V_{k-1}<0$ and then
$\Delta F_{k-1}<0$.

From (\ref{0k-1}) and $0\le \Delta\tilde Q_k \le A_o\Delta\tilde V_k$
we have
\begin{eqnarray}
0<\Delta\tilde F_k \le (1+KA_o) \piu{\Delta\tilde V_k} \le \lambda
\biggl( \meno{\Delta V_{k-1}} - \sum_{\ell = 1}^{k-2}\piu{\Delta
V_\ell} \biggr)\,.
\label{eq:ll}
\end{eqnarray}
We prove now that
\begin{eqnarray}
\meno{\Delta Q_{k-1}} - \sum_{\ell = 1}^{k-2}\piu{\Delta Q_\ell}
\ge 0\,.
\label{eq:mm}
\end{eqnarray}
We have only to consider the case in which $[\Delta Q_\ell]_+>0$ for
some $\ell\le k-2$; but in this case $[\Delta Q_{k-1}]_- - \sum_{\ell
  = 1}^{k-2}[\Delta Q_\ell]_+ = L_{cd}^+ (|\alpha_3| + |\beta_3| -
|\eps_3|)\ge 0$ because of (\ref{eq:333}).

Therefore (\ref{h=k-1b-XXX}) for $\tau\in I_{k-1}$ follows from
(\ref{eq:ll}) and (\ref{eq:mm}).

\medskip

Second, assume $\tau\in J_{k-1}$; we prove that
\begin{eqnarray}\label{h=k-1b-XXX-caso-J}
\piu{\Delta\tilde F_k} \le \mu
\meno{\Delta F_{k-1}} \,.
\end{eqnarray}
Indeed, if the reflected wave is a physical wave then, under the
notation of Proposition~\ref{Delta-F-2wave},
$$ \Delta V_{k-1}
\le \xi \frac{|\delta_1\delta_2|}{2}\,,\qquad \Delta Q_{k-1} \le -
|\delta_1\delta_2| + \frac{|\delta_1\delta_2|}{2} A_o = -
\frac{|\delta_1\delta_2|}{2} (2-A_o)
$$
so that $\Delta F_{k-1}\le \left[\xi - K(2-A_o)\right]
|\delta_1\delta_2|/2 <0$ because of (\ref{cond_su_chi_xi-NEW}).
Then (\ref{h=k-1b-XXX-caso-J}) follows since $\Delta\tilde F_k =
\Delta F_k >0$ and
\[
\piu{\Delta\tilde F_k}
\le \frac{|\delta_1\delta_2|}{2} (\xi+KA_o)
= \frac{|\delta_1\delta_2|}{2}\left[K(2-A_o) - \xi\right] \cdot \lambda_2
\le \lambda_2 \meno{\Delta F_{k-1}} \le \mu \meno{\Delta F_{k-1}}\,.
\]
If the reflected wave is a non-physical wave we have, under the
notation of Proposition~\ref{lem:np},
\begin{eqnarray*}
0< \Delta F_k = \Delta V_k \le \knp C_o|\delta\delta_2|\,,\qquad
\Delta V_{k-1}=0\,,\qquad \Delta Q_{k-1} =  - |\delta\delta_2|
\end{eqnarray*}
and then
\begin{eqnarray*}
\meno{\Delta F_{k-1}} = K|\delta\delta_2|\,,\qquad \piu{\Delta F_k}\le
\frac{\knp C_o}{K}\meno{\Delta F_{k-1}}\le \mu\meno{\Delta F_{k-1}}\,.
\end{eqnarray*}
The estimate (\ref{h=k-1b-XXX-caso-J}) is then completely proved.
From (\ref{h=k-1b-XXX-caso-J}) we get (\ref{h=k-1b-XXX}) since no
waves of order $\le k-2$ are involved in the interaction.

\item Finally, let $h\ge k$.  We first consider $\tau\in I_{h}$,
$h\ge k$. If $\min\{k_\alpha, k_\beta\}\ge k$, then $\Delta \tilde
V_k =\Delta L_\xi < 0$. If $\min\{k_\alpha, k_\beta\}\le k-1$,
assume $k_\alpha \ge k$ and $k_\beta \le k-1$; then $\Delta \tilde
V_k\leq \xi|\eps_1| - |\alpha_3|\leq (\xi d -1) |\alpha_3|< 0$ by
(\ref{cond_su_chi_xi})${}_1$. This proves $\Delta \tilde V_k < 0$
and then $\Delta \tilde F_k <0$.

From (\ref{Delta_W_piu_V}) we deduce that $\sum_{\ell=1}^{k-1}\Delta
F_{\ell} < [\Delta \tilde F_k]_- $. Since at most one non-zero term is
present in the first sum, (\ref{Delta_tilde_V_k-meno}) follows. In the
case $h=k$ we have $\Delta \tilde V_k = [\Delta V_{k+1}]_+ + \Delta
V_k < 0$; hence $\Delta V_k < 0$ and then $\Delta F_k<0$.

Now assume $\tau\in J_{h}$, $h\ge k$. Then $ \Delta \tilde F_k =
\Delta F <0$ and no waves of order $<k$ are present, so
(\ref{Delta_tilde_V_k-meno}) holds. If $h=k$ and the reflected wave is
physical, from the proof of Proposition~\ref{Delta-F-2wave} and
(\ref{cond_su_chi_xi-NEW}) we find that
$$
\Delta F_k = \Delta V_k + K  \Delta Q_k \le
\frac{|\delta_1\delta_2|}2 \left(\xi - 2K + KA_o\right)<0\,.
$$
If $h=k$ and the reflected wave is non-physical then $\Delta V_k =
0$, $\Delta Q_k <0$ and $\Delta F_k<0$.
\end{remunerate}
\end{proof}
\noindent Summarizing, for $\tau\in {\cal T}$ we have the
following table:
\[
\begin{tabular}{r|cccc}
${\cal T}_h$; $h$ & $\le k-2$ & $k-1$ & $k$ & $\ge k+1 $
\\
\hline $\Delta F_k$ & $0$ & $+$ & $-$ & $\pm$
\\
$\Delta \tilde F_k$ & $0$ & $+$ & $-$ & $-$
\end{tabular}
\]
\noindent We write $\tilde F_k^\pm(t) = \sum_{\tau \le t}[\Delta
\tilde F_k(\tau)]_\pm$ for $k\ge 2$. For simplicity the time $\tau$ in
such sums is omitted.

%

\begin{lemma}
  Under the assumptions of Proposition~\ref{lem:01} we have
\begin{eqnarray}
\tilde F_2^+(t) & \le & \mu\left(L_\xi(0) + K Q(0)\right) +
\sum_{{\cal T}_{h},\,h\ge 2}\piu{\Delta F_1}\label{03}
\\
\tilde F_k^+(t) & \le & \mu \Bigl(\tilde F_{k-1}^+(t) +
\sum_{{\cal T}_{h},\,h\ge k}\piu{\Delta F_{k-1}} - \sum_{{\cal
T}_{k-1}}\sum_{\ell = 1}^{k-2}\piu{\Delta F_\ell}\Bigr)\,, \qquad
k\ge 3\,. \label{04}
\end{eqnarray}
\end{lemma}

\begin{proof}
  By Proposition~\ref{lem:01}, see also the table above, the
  functional $F_{k}$ increases at times $\tau \in {\cal T}_{k-1}$,
  decreases at times $\tau \in {\cal T}_{k}$, while it does not have a
  given sign at times $\tau \in {\cal T}_{h}$, with $h\ge k+1$.

First, by summing up (\ref{h=k-1b-XXX}) we obtain
\begin{eqnarray}
\tilde F_k^+(t)
& \le &
\mu \sum_{{\cal T}_{k-1}} \Bigl(
\meno{\Delta F_{k-1}} - \sum_{\ell = 1}^{k-2}\piu{\Delta F_\ell}\Bigr)
\label{0cor0}
\end{eqnarray}
for $k\ge 2$, where the last term in (\ref{0cor0}) is missing if
$k=2$.

Recall now that $F_1(0) = L_\xi(0) + K Q(0)$; therefore
\begin{eqnarray*}
F_1(t) &\leq & L_\xi(0) + K Q(0) - \sum_{{\cal T}_1} \meno{\Delta
F_{1}} + \sum_{{\cal T}_h,\,h\ge{2}} \piu{\Delta F_{1}}
\end{eqnarray*}
and then
\begin{eqnarray} \sum_{{\cal T}_1}\meno{\Delta F_{1}}
&\leq & L_\xi(0) + K Q(0) + \sum_{{\cal T}_h,\, h\geq 2} \,
\piu{\Delta F_{1}}\,. \label{Delta_V_1-meno-NEW}
\end{eqnarray}
On the other hand $F_k(0) = 0$ for $k\ge 2$; from
Proposition~\ref{lem:01} we have
\begin{eqnarray*}
F_k(t)
&\leq &
\sum_{{\cal T}_{k-1}} \piu{\Delta F_{k}} - \sum_{{\cal T}_{k}} \meno{\Delta F_{k}}
+ \sum_{{\cal T}_h,\, h\geq k+1} \piu{\Delta F_{k}}\,.
\end{eqnarray*}
Moreover
$$
\sum_{{\cal T}_{k-1}} \piu{\Delta F_{k}}
=
\sum_{{\cal T}_{k-1}} \piu{\Delta \tilde F_{k}}= \tilde F_{k}^+(t)
$$
and then
\begin{eqnarray}
\sum_{{\cal T}_{k}}\meno{\Delta F_{k}} & \leq & \tilde F_{k}^+(t)
+ \sum_{{\cal T}_h,\, h\geq k+1}\, \piu{\Delta F_{k}}\,.
\label{Delta_V_k-meno-NEW}
\end{eqnarray}
From (\ref{0cor0}), (\ref{Delta_V_1-meno-NEW}),
(\ref{Delta_V_k-meno-NEW}) we get (\ref{03}), (\ref{04}).
\end{proof}

\smallskip

\begin{proposition}[A contraction property]\label{prop:piubella}
  Under the assumptions of Proposition~\ref{lem:01}, for any $t\ge0 $
  and $k\ge 1$ we have
\begin{eqnarray}
\tilde V_k(t)\le \tilde F_k(t) & \le & \mu^{k-1}\cdot \left( L_\xi(0) + K Q(0)\right)\,. \label{piu-bella}
\end{eqnarray}
\end{proposition}
\begin{proof}
  The estimate (\ref{piu-bella}) holds for $k=1$ because $\tilde
  F_1(t) = F(t) \le L_\xi(0) + K Q(0)$. Next we prove by induction on
  $k\ge 2$ that for any $t$
\begin{eqnarray}
\tilde F_k^+(t)
& \le &
\mu^{k-1} \left( L_\xi(0) + K Q(0)\right)
+ \sum_{{\cal T}_{h},\, h\ge k}\sum_{\ell = 1}^{k-1}\piu{\Delta F_\ell}\,.
\label{prima-induzione}
\end{eqnarray}
Since by summing up (\ref{Delta_tilde_V_k-meno}) we obtain
\begin{eqnarray}
\tilde F_k^-(t) & \ge & \sum_{{\cal T}_{h},\,h\ge
k}\sum_{\ell = 1}^{k-1}\piu{\Delta F_\ell} \label{0cor1}
\end{eqnarray}
then (\ref{piu-bella}) will follow from (\ref{prima-induzione})
for any $k\ge2$ because of (\ref{0cor1}).

Formula (\ref{prima-induzione}) for $k=2$ reduces to (\ref{03}).
Next, assume that (\ref{prima-induzione}) holds for some $k\ge 2$.  By
(\ref{04}) and the induction assumption
\begin{eqnarray*}
\lefteqn{\tilde F_{k+1}^+(t) ~ \le ~ \mu\Bigl( \tilde
F_{k}^+(t) + \sum_{{\cal T}_{h},\, h\geq k+1} \, \piu{\Delta F_{k}}
- \sum_{{\cal T}_{k}}\sum_{\ell = 1}^{k-1}\piu{\Delta
F_{\ell}}\Bigr)}
\\
& \le &
\mu^k \left( L_\xi(0) + K Q(0)\right) + \mu\Bigl( \sum_{{\cal T}_{h},\, h\ge k}\sum_{\ell = 1}^{k-1}\piu{\Delta
F_\ell} + \sum_{{\cal T}_{h},\, h\geq k+1} \piu{\Delta F_{k}}
- \sum_{{\cal T}_{k}}\sum_{\ell = 1}^{k-1}\piu{\Delta F_{\ell}}
\Bigr)
\\
& \le &
{\mu^k} \left( L_\xi(0) + K Q(0)\right) + \mu \sum_{{\cal T}_{h},\,
h\ge k+1}\sum_{\ell = 1}^{k} \piu{\Delta F_\ell}\,.
\end{eqnarray*}
Since $\mu<1$, we get (\ref{prima-induzione}) for $k+1$.
\end{proof}

\smallskip

\begin{remark}\label{rem:L2=0}\rm
  We comment now the case $A_o=0$. In this case system
  (\ref{eq:system}) reduces to the $p$-system with pressure law given
  by (\ref{eq:pressure}), for fixed $\lambda$. According to our
  front-tracking algorithm, stationary and non-physical waves do not
  appear and so $L_{cd}=Q=0$; the algorithm reduces to the one introduced
  in \cite{AmadoriGuerra01}.
Then, Proposition \ref{lem:01} holds with $\tilde
  V_k$ and $1/\xi$ replacing $\tilde F_k$ and $\mu$, respectively, and
  at last (\ref{piu-bella}) reads
\begin{eqnarray*}
\tilde V_k(t) & \le & \frac{1}{\xi^{k-1}} \cdot L_\xi (0)\,.
\end{eqnarray*}
\end{remark}

\smallskip

Next, we conclude that the total strength of all non-physical
waves is small by proceeding as in \cite[page 142]{Bressanbook}.
Recall the notation in Section \ref{sec:app_sol}. First, using
(\ref{eq:opt}), we find that the sequence $\{v^\nu\}$ related to
the initial data $(v_o^\nu,u_o^\nu,\lambda_o^\nu)$ is uniformly
bounded and in particular uniformly bounded away from $0$; then
the eigenvalues $e_1$ and $e_3$ are bounded and this makes
possible the choice of a suitable $\hat s$. We have two more
parameters $\eta$, $\rho$ to be chosen. Fix $\eta>0$ with the
condition $\eta=\eta_\nu\to 0$ as $\nu\to\infty$ and estimate the
total number of waves of order $<k$. We have
\begin{eqnarray*}
\sum_{\gamma \in {\cal N\!P}} |\gamma| (t) &\le&
\tilde V_k(t) ~+~ \sum_{\gamma \in {\cal N\!P},\,k_\gamma<k} |\gamma|
(t)\\
&\le&
\mu^{k-1}\cdot \left( L_\xi(0) + K Q(0)\right)  ~+~ C_o\rho \cdot
[\mbox{number of fronts of order $<k$}]
~\le~ \frac{1}{\nu}
\end{eqnarray*}
by choosing $k$ sufficiently large to have the first term $\le
1/(2\nu)$ and then choose $\rho$ small enough to have the second term
$\le 1/(2\nu)$.

\bigskip

We accomplish now the proof of Theorem \ref{thm:main}. Define
first
\begin{eqnarray}
k(m) & = & \frac{1-\sqrt{d(m)}}{2-\sqrt{d(m)}}\,. \label{eq:kM}
\end{eqnarray}
From the properties of the function $d(m)$ stated in Remark
\ref{rem:d(m)} we see that $k(0)=1/2$ and that $k(m)$ is decreasing,
tending to $0$ for $m\to+\infty$. The assumption (\ref{hyp2}) implies
that (\ref{bound_su_A_o-NEW}) holds.

Now, by hypotheses (\ref{hyp1}) it follows that we can choose $\xi$
such that
$$
\frac{1}{2}\tv \log(p_o) + \frac{1}{2\inf a_o}\tv(u_o) ~<~ \frac
m{2\xi-1} ~<~ \bigl(1 - 2A_o\bigr)m
$$
and that (\ref{cond_su_chi_xi-NEW})${}_1$ holds.
Hence, using (\ref{stima-onde-PdR}) and \emph{(i)} in Section
\ref{sec:app_sol}, we have
\begin{eqnarray}
\nonumber
L(0+)
& \le&  \frac{1}{2}\tv \log(p_o) + \frac{1}{2\inf a_o}\tv(u_o) < \frac
m{2\xi-1}
\end{eqnarray}
so that the hypotheses (\ref{eq:super_cond}) of
Proposition~\ref{prop:very_last} holds.  Theorem \ref{thm:main}
now follows along the lines of \cite[\S 7.4]{Bressanbook}.


\appendix

\section{The weighted total variation}\label{app:WTV}
In this Appendix we prove Proposition \ref{prop:wtv}. Remark that
the map $d(a,b) \doteq \frac{|a-b|}{a+b}$ is a distance on
$\reali_+$, as one can easily prove.

\smallskip

We start with the proof of the inequality on the right in
(\ref{eq:smaller}). It is enough to prove that
\begin{equation}\label{eq:greater}
2\sum_{j=1}^n\frac{|f(x_j)-f(x_{j-1})|}{f(x_j)+f(x_{j-1})}\le
\sum_{j=1}^n\left|\log f(x_j)-\log f(x_{j-1})\right|\,.
\end{equation}
We claim that
\begin{equation}\label{eq:newlog}
\log t\ge \frac{2(t-1)}{t+1}\quad\hbox{ for }t\ge1
\end{equation}
where the inequality is strict if $t>1$. To prove the claim it is
sufficient to notice that the function $\phi(t) = \log
t-2\frac{t-1}{t+1}$ vanishes in $1$ and $\phi'(t) =
\frac{(t-1)^2}{t(t+1)^2}>0$ if $t>1$.

We apply (\ref{eq:newlog}) to $t=x/y$ for $0<y\le x$ and arguing
by symmetry deduce
\[
|\log x-\log y|\ge \frac{2|x-y|}{x+y}\,, \quad\hbox{ for every }
x,y > 0\,.
\]
Then (\ref{eq:greater}) follows. The proof of the inequality on
the left in (\ref{eq:smaller}) is analogous, starting from the
inequality
\begin{equation}\label{eq:newnewlog}
\frac{1}{t}\log t\le \frac{2(t-1)}{t+1}\quad\hbox{ for }t\ge1
\end{equation}
with strict inequality if $t>1$.

Assume now that $f\in C(\reali)$. To show that $\wtv(f) =
\tv\left(\log(f)\right)$, we have to prove that the inequality
\begin{equation}\label{eq:converse}
\tv\left(\log(f)\right)\le
2\sup\sum_{j=1}^n\frac{|f(x_j)-f(x_{j-1})|}{f(x_j)+f(x_{j-1})}
\end{equation}
holds for any interval $[a,b]\subset \reali$. Consider any
partition $\{x_o,x_1,\ldots,x_n\}$ of the interval $[a,b]$. By the
mean value theorem for some $\zeta_j$ between $f(x_j)$ and
$f(x_{j-1})$, and by the intermediate value theorem applied to
$f$, we get
\begin{eqnarray*}
\left|\log\left(f(x_j)\right) - \log\left(f(x_{j-1})\right)\right|
& = & \frac{\left|f(x_j)-f(x_{j-1})\right|}{\zeta_j}
\\
& = & \frac{f(x_j)+f(x_{j-1})}{2 f(\eta_j)} \cdot
2\frac{\left|f(x_j)-f(x_{j-1})\right|}{f(x_j)+f(x_{j-1})}
\end{eqnarray*}
for some $\eta_j\in[x_{j-1},x_j]$. We exploit again the continuity
of $f$ in $[a,b]$. On one hand its image is compact, then
$\min_{[a,b]}f=m > 0$. On the other $f$ is uniformly continuous in
$[a,b]$, so that for any $\epsilon>0$ there exists
$\delta_\epsilon > 0$ such that $|f(x)-f(y)|<\epsilon$ if $|x-y| <
\delta_\epsilon$, for $x,y\in[a,b]$.

Fix now any $\epsilon>0$; without loss of generality we can
consider partitions of the interval $[a,b]$ of mesh less than
$\delta_\epsilon$. Assume for instance $f(x_{j-1})\le f(x_j)$;
then from the inequalities
\begin{eqnarray*}
f(x_j)-f(x_{j-1}) ~<~\epsilon\,,\qquad f(x_{j-1})\le  f(\eta_j)
\le f(x_j)
\end{eqnarray*}
it follows
\[
\frac{f(x_j)+f(x_{j-1})}{2 f(\eta_j)} ~\le~
\frac{2f(x_{j-1})+\epsilon}{2f(x_{j-1})}~\le~
1+\frac{\epsilon}{m}\,.
\]
The inequality (\ref{eq:converse}) follows by remarking that then
for any partition of mesh less than $\delta_\epsilon$
\begin{eqnarray*}
\sum_{j=1}^n\left|\log\left(f(x_j)\right) -
\log\left(f(x_{j-1})\right)\right| & \le &
\left(1+\frac{\epsilon}{m}\right)
2\,\sum_{j=1}^n\frac{|f(x_j)-f(x_{j-1})|}{f(x_j)+f(x_{j-1})}\,.
\end{eqnarray*}

The proof of Proposition~\ref{prop:wtv} is complete.

\smallskip

\begin{remark}\label{rem:log-f}\rm
  Observe that, if $\tv\left(\log(f)\right)<\infty$ and $f$ is
  discontinuous, then the inequality on the right in (\ref{eq:smaller}) is strict,
  because of the strict inequality in (\ref{eq:newlog}). For example,
if $f$ has a single jump and assumes the values $c>0$ and $d>0$,
then $\wtv(f) = 2\frac{|c-d|}{c+d} < |\log c-\log d| =
\tv(\log(f))$.

\smallskip
Remark moreover that $\wtv$ and $\tv$ are not equivalent, in the
sense that there does not exist a positive constant $C$ such that
$C \cdot \tv\left(\log(f)\right)\le \wtv(f)$. This follows from
the fact that clearly the inequality $C\log t\le
\frac{2(t-1)}{t+1}$ does not hold for every $t\ge1$.
\end{remark}


\section{Shock-rarefaction interactions}\label{app:SS}
In this Appendix we consider a particular case of Lemma
\ref{lem:shock-riflesso}. Actually, this is the only case needed in
order to define a decreasing functional,
Section~\ref{sec:interactions}; however, we needed further analysis
for the control and treatment of the non-physical waves.

\smallskip

\begin{lemma}[The case $SR,\, RS\to
SS$]\label{lem:second} Consider the interaction of a shock
$\alpha_i$ and a rarefaction $\beta_i$ of the same family,
$i=1,3$, producing two outgoing shocks $\eps_1$, $\eps_3$. Then
there exists a smooth function $B$ satisfying $|\alpha_i| \le
B(\alpha_i)\le \min\{\sinh(|\alpha_i|),2|\alpha_i|\}$ such that
\begin{equation}\label{eq:boundforb}
0<\beta_i \le B(\alpha_i)\,.
\end{equation}
Moreover, assume
\begin{eqnarray}
|\alpha_i| & \le & m
\label{alpham}
\end{eqnarray}
for some $m>0$ and denote $c = c(m) =
\frac{\cosh(m)-1}{\cosh(m)+1}$. Then both the variation of shock
waves and the reflected wave $\eps_j$, $j\ne i$, $j=1,3$, are
estimated by the interacting rarefaction as
\begin{eqnarray}
|\eps_1| + |\eps_3| - |\alpha_i| & \le & (2c-1)\cdot |\beta_i|\,,
\label{eq:var_shocks}
\\
|\eps_j| & \le & c \cdot|\beta_i|\,.\label{eq:eps1c}
\end{eqnarray}
\end{lemma}

\begin{proof}
  We focus on the case $i=3$, $j=1$, see Figure \ref{fig:inter3133}(b).
  Therefore we consider $\alpha_3<0$, $\beta_3>0$ and $\eps_1<0$,
  $\eps_3<0$.  Then (\ref{tre-uno}), (\ref{tre-due}) become
\begin{eqnarray}
|\eps_1|-|\eps_3| & = & |\beta_3|-|\alpha_3|\label{eq:IVa}
\\
\sinh(|\eps_1|)+ \sinh(|\eps_3|) & = & \sinh(|\alpha_3|) -
|\beta_3| \,.
\label{eq:IVb}
\end{eqnarray}
From the second equation $|\eps_1|< |\alpha_3|$, $|\eps_3|<
|\alpha_3|$ and $|\beta_3|< \sinh(|\alpha_3|)$; using the first
equation and $|\eps_3|< |\alpha_3|$, we get $|\eps_1|<|\beta_3|$.
Therefore in conclusion
\begin{eqnarray}\label{inequality-1}
|\eps_1|< \min\{|\alpha_3|,|\beta_3|\} ,\qquad\ |\eps_3|< |\alpha_3|,\qquad
|\beta_3|< \sinh(|\alpha_3|)\,.
\end{eqnarray}

\smallskip
\noindent {\bf Step 1: notation.} We set
$x=|\beta_3|$, $y=|\eps_1|$, $z=|\alpha_3|$, so that
\begin{equation}\label{eq:eps3pos}
|\eps_3|=y-x+z\,.
\end{equation}
Under this notation, (\ref{eq:IVb}) writes as
\begin{equation}\label{eq:impl_eq}
F(x,y;z) = \sinh y+ \sinh(y-x+z) - \sinh z + x = 0
\end{equation}
for $x\ge0$, $y\ge0$, $z\ge0$, $y-x+z\ge0$. By
(\ref{inequality-1}), any solution of (\ref{eq:impl_eq}) satisfies
\begin{equation}\label{eq:prop_y}
y<z\,,\qquad y<x\,,\qquad x<\sinh(z)\,.
\end{equation}


\begin{figure}[htbp]
\begin{picture}(100,60)(-20,-30)
\setlength{\unitlength}{1pt}
\put(160,0){
\put(0,0){\line(-1,1){30}}\put(-50,30){\makebox(0,0){$y=|\eps_1|$}}
\put(0,0){\line(1,1){30}}\put(70,30){\makebox(0,0){$y-x+z=|\eps_3|$}}
\put(0,0){\line(-2,-1){40}}\put(-60,-20){\makebox(0,0){$z=|\alpha_3|$}}
\put(0,0){\line(-2,-3){20}}\put(5,-25){\makebox(0,0){$x=|\beta_3|$}}
}
\end{picture}
\caption{\label{fig:inter_xyz}{Interactions.}}
\end{figure}


\smallskip

\noindent {\bf Step 2: the threshold.} Observe that, despite the
last inequality in (\ref{inequality-1}), we may well have
$|\beta_3|>|\alpha_3|$, that is, $x>z$. Consider in fact the limit
case of $\eps_3=0$: we have $y  =  x-z>0$ and $\sinh(y)  =
\sinh(z) -x$ that give
\begin{eqnarray*}
\sinh(x-z) & = & \sinh(z) -x,\qquad x> z \,.
\end{eqnarray*}
The last equality is the relation needed for $\beta_3$, $\alpha_3$
in order to have that the shock and rarefaction cancel out
exactly, giving rise only to a wave of the opposite family.
Observe that the size of the rarefaction must be larger than the
one of the shock. Under the notations above, the threshold curve
separating the case of the outgoing waves $S_1S_3$ from the case
$S_1R_3$ is given by
\begin{eqnarray}\label{soglia}
f(x,z) = \sinh(x-z) - \sinh(z)+ x = 0\,.
\end{eqnarray}
Since $f_x=\cosh(x-z)+1>0$ and $f(z,z)<0$, the implicit equation
$f(x,z)=0$ is solved by $x=x_o(z)\ge z$ with
$x_o'(z)=\frac{\cosh(x-z)+\cosh(z)}{\cosh(x-z)+1}>0$ for every
$z\ge0$; the curve has for tangent at $(0,0)$ the line $z=x$.
Observe that $x_o(z) \leq 2z$ because $f(2z,z)=z\geq 0$. In
conclusion
\begin{equation}\label{eq:xoz}
z\le x_o(z) \leq 2z\,.
\end{equation}
This estimate and the third inequality in (\ref{inequality-1})
prove (\ref{eq:boundforb}) for $B(\alpha_i)=x_o(|\alpha_i|)$. We
can prove more than (\ref{eq:xoz}), that is,
\begin{equation}\label{eq:xoz-2z}
\lim_{z\to+\infty}\left(x_o(z) - 2z\right)=0\,.
\end{equation}
Indeed we show that the inequality $x_o(z)> 2z - q$ holds for large $z$ and $q>0$ . This
follows from $ f(2z-q,z) = \sinh(z-q) - \sinh z + 2z -q \sim e^z
\left(e^{-q} -1\right)/2 \to -\infty $ for $z\to\infty$.

Remark that
the fact that $\eps_3$ is a shock
implies that
\begin{eqnarray}\label{soglia-2}
f(x,z) = \sinh(x-z) - \sinh(z)+ x < 0\,.
\end{eqnarray}


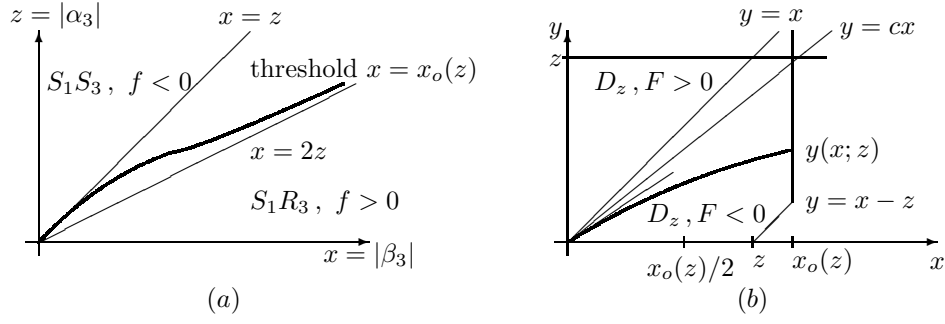
\begin{figure}[t]
\begin{picture}(120,120)(-20,-25)
\setlength{\unitlength}{1pt}

\put(0,0){
\put(0,0){\vector(1,0){125}}
\put(125,-10){\makebox(0,0)[b]{$x=|\beta_3|$}}
\put(0,0){\line(-1,0){5}} \put(0,0){\vector(0,1){80}}
\put(-11,86){\makebox(0,0)[l]{$z=|\alpha_3|$}}
\put(0,0){\line(0,-1){5}}
\put(0,0){\line(1,1){80}}
\put(68,85){\makebox(0,0)[l]{$x=z$}}
\put(0,0){\line(2,1){120}} \put(80,35){\makebox(0,0)[l]{$x=2z$}}
\put(0,0){\thicklines{\qbezier(0,0)(22,24)(50,34)}}
\put(0,0){\thicklines{\qbezier(50,34)(60,35)(115,60)}}
\put(80,65){\makebox(0,0)[l]{threshold $x=x_o(z)$}}
\put(3,60){\makebox(0,0)[l]{$S_1S_3\,,\ f<0$}}
\put(80,15){\makebox(0,0)[l]{$S_1R_3\,,\ f>0$}} }
\put(70,-22){\makebox(0,0){$(a)$}}


\put(200,0){ \put(0,0){\vector(1,0){140}}
\put(140,-10){\makebox(0,0)[b]{$x$}} \put(0,0){\line(-1,0){5}}
\put(0,0){\vector(0,1){80}} \put(-7,80){\makebox(0,0)[l]{$y$}}
\put(0,0){\line(0,-1){5}}

\put(0,0){\line(1,1){80}} \put(65,85){\makebox(0,0)[l]{$y=x$}}
\put(70,0){\line(1,1){15}} \put(90,15){\makebox(0,0)[l]{$y=x-z$}}
\put(70,-2){\line(0,1){4}}
\put(70,-7){\makebox(0,0)[l]{$z$}}
\put(-2,70){\line(1,0){100}}
\put(-7,69){\makebox(0,0)[l]{$z$}}
\put(85,15){\line(0,1){65}}
\put(85,-7){\makebox(0,0)[l]{$x_o(z)$}} \put(44,-2){\line(0,1){4}}
\put(30,-10){\makebox(0,0)[l]{$x_o(z)/2$}}
\put(85,-2){\line(0,1){4}}
\put(10,60){\makebox(0,0)[l]{$D_z\,, F>0$}}
\put(30,10){\makebox(0,0)[l]{$D_z\,, F<0$}}
\put(0,0){\line(3,2){40}} \put(0,0){\line(5,4){100}}
\put(103,80){\makebox(0,0)[l]{$y=cx$}}
\put(0,0){\thicklines{\qbezier(0,0)(40,25)(85,35)}}
\put(90,35){\makebox(0,0)[l]{$y(x;z)$}}
\put(70,-22){\makebox(0,0){$(b)$}}
}


\end{picture}

\caption{\label{fig:lax}{\emph{(a)}, the threshold curve $f(x,z)=\sinh(x-z)  -
\sinh(z) +x=0$; \emph{(b)}, the domain $D_z$ and the function
$y=y(x;z)$.}
}
\end{figure}


\noindent {\bf Step 3: the amount of shocks can increase.} From (\ref{eq:IVa}) we have
\begin{equation}\label{eq:equals}
|\eps_1|+|\eps_3|-|\alpha_3| = 2|\eps_1|-|\beta_3|\,.
\end{equation}
We prove now that the inequality $|\eps_1|+|\eps_3|-|\alpha_3| <
0$, or equivalently $|\eps_1| < \frac12 |\beta_3|$, does not hold
if $m$ is large.

The equation giving $|\eps_1|=y$ in terms of $|\beta_3|=x$, for a
given parameter $|\alpha_3|=z$, is (\ref{eq:impl_eq}), to be
considered in the domain
\[
D_z=\left\{(x,y)\colon\  0\leq x \leq x_o(z) \,,\
y\ge\max\{0,x-z\} \right\}\,, \qquad \hbox{ for } z\ge0
\]
where $x_o(z)$ satisfies the equation (\ref{soglia}). Since $F_y =
\cosh y + \cosh(y-x+z)>0$, the implicit equation
(\ref{eq:impl_eq}) defines a function $y=y(x;z)$ with $0\le
y(x;z)\le x$ and $y(x;z)\le z$, see (\ref{eq:prop_y}). Remark that
$F(x,x;z) = \sinh x + x >0$. Moreover $F(x,0;z) = \sinh(z-x) -
\sinh z + x$ so that $F(0,0;z) = 0$; by $F_x=1-\cosh(z-x)<0$, we
deduce $F(x,0;z) <0$ if $x>0$. By the implicit function theorem we
have $y'(x;z)\ge0$,
\begin{equation}\label{eq:y'}
y'(0;z) = \frac{\cosh(z)-1}{\cosh(z)+1}\in[0,1)
\end{equation}
and $y''(0;z)=-\frac{4\sinh(z)}{(1+\cosh(z))^2}\le0$. The function
$\frac{\cosh(z)-1}{\cosh(z)+1}$ is increasing and then for
$z\in[0, m]$ its maximum is $\frac{\cosh(m)-1}{\cosh(m)+1}$; this
quantity is strictly larger than $1/2$ if $ m > \log(3+2\sqrt2)$.
Thus in general the estimate $y(x;z) < x/2$ cannot hold.

\smallskip

\noindent {\bf Step 4: proof of the estimate.} From
(\ref{eq:equals}) we see that $|\eps_1| + |\eps_3| - |\alpha_3|
\le (2c-1)\cdot |\beta_3| \iff |\eps_1| \le c \cdot|\beta_3|$,
that is, (\ref{eq:var_shocks}) and (\ref{eq:eps1c}) are
equivalent; we shall prove (\ref{eq:eps1c}). To bypass the study
of the function $y(x;z)$ we define
\begin{eqnarray*}
\Phi(x;z,c)=F\left(x, c x;z\right) &=& \sinh (cx)+
\sinh\left(z-(1-c)x\right) - \sinh z + x\,.
\end{eqnarray*}
If $1/2\leq c<1$ then $z>(1-c)x$ and $\Phi(0;z,c)=0$\,,
\begin{eqnarray*}
\Phi_x(x;z,c) & = & 1 + c \cosh(cx) - (1-c) \cosh\left(z -(1-c)
x\right)\,,
\\
\Phi_{xx}(x;z,c) & = & c^2 \sinh(cx) +  (1-c)^2 \sinh\left(z
-(1-c) x\right) >0\,.
\end{eqnarray*}
Therefore the function $x\to \Phi_x(x;z,c)$ is increasing and then
$\Phi(x;z,c)\ge0$ if $$\Phi_x(0;z,c) = (1 + c) - (1-c)
\cosh(z)>0\,,$$ that is if $\cosh(z) \le \frac {1 + c}{1-c}$; this
is just (\ref{alpham}). Then $y(x;z) \le cx$ for all
$x\in\left(0,x_o(z)\right)$ and so (\ref{eq:eps1c}) is proved.
\end{proof}

\smallskip

\begin{remark}\label{rem:log}\rm
From (\ref{eq:y'}) we see that condition (\ref{alpham}) is
equivalent to the geometric condition $y'(0;z) = (\cosh
z-1)/(\cosh z+1)<c$. Moreover, as we noticed in the above proof,
condition (\ref{alpham}) is equivalent to
\begin{equation}\label{eq:cosh1c}
\cosh(|\alpha_i|) \le \frac{1+c}{1-c}\,,
\end{equation}
which, in turn, is equivalent to
\begin{equation}\label{eq:thalpha}
|\alpha_i| \leq \log \frac{(1+\sqrt{c})^2}{1-c}\,.
\end{equation}
From the definition of the strength, one has that
$|\alpha_i|=(1/2) \log(v_{max}/v_{min})$, where
$v_{max}=\max\{v_\ell,v_r\}$, $v_{min}=\min\{v_\ell,v_r\}$, being
$v_\ell$, $v_r$, respectively, the left and right values of $v$
for the wave of size $\alpha_i$. Hence (\ref{eq:thalpha}) is
equivalent to
\begin{eqnarray*}
\sqrt{\frac{v_{max}}{v_{min}}} \leq \frac {\left(1 +
\sqrt{c}\right)^2}{1-c}\,.
\end{eqnarray*}
\end{remark}

\smallskip

\begin{remark}\rm
  In the proof above we showed that $\Delta L_{{\rm shocks}} =
  |\eps_1|+|\eps_3|-|\alpha_3| = 2|\eps_1|-|\beta_3|$ may be positive,
  differently from Nishida's paper, where it is always decreasing.
  This depends on the definition of the wave strengths, which was
  imposed to us in order to have good estimates when dealing with
  interactions with the $2$ waves. In any case $\Delta L =
  \Delta L_{{\rm shocks}} + \Delta L_{{\rm rarefactions}}
  \le0$.
\end{remark}

\smallskip
\begin{remark}\rm
  Under the notation and assumptions of Lemma~\ref{lem:second} we
  verify that
\begin{equation}\label{eq:eps1_c_alpha3}
|\eps_j|\leq  c \cdot |\alpha_i|\,.
\end{equation}
This estimate, together with (\ref{eq:eps1c}), allows us to obtain (in
a special case) the analogue of (\ref{eq:chi_def}) with $c$ in place
of $d$.

The proof makes use of a numerical computation. As in that lemma we
consider the case $j=1$, $i=3$.  Let $z>0$ be fixed and $0\leq x\leq
x_o(z)$. From the proof of Lemma \ref{lem:second} we deduce
$y=y(x,z)$; for $z$ fixed the function $x\to y(x,z)$ is increasing.
Define then
\[
Y(z) = y(x_o(z),z) = x_o(z)-z\,.
\]
In order to prove (\ref{eq:eps1_c_alpha3}) it is sufficient to
prove that
\[
Y(z)\leq cz \quad \hbox{if} \quad \cosh(z)\le \frac{1+c}{1-c}\,.
\]
Remark that from (\ref{eq:xoz-2z}) we know that $Y(z)/z \to 1$ for
$z\to+\infty$; however the constraint (\ref{eq:cosh1c}) implies
that $z$ is bounded. The inequality $Y(z)\leq cz$ is equivalent to
$x_o(z)-z \le cz$, i.e., $x_o(z)\le (1+c)z$. Therefore we need to
prove that
\begin{equation}\label{eq:fabrizio}
\phi(z,c) \doteq \sinh(cz) - \sinh(z)+(1+c) z \ge0 \quad \hbox{ if
}\quad \cosh(z)\le \frac{1+c}{1-c}\,.
\end{equation}
Notice that $\cosh(z)\le \frac{1+c}{1-c}$ means $0\le z\le z_c$
for $z_c = \log(\frac{1+c}{1-c} +\sqrt{(\frac{1+c}{1-c})^2-1}\,)$.
Formula (\ref{eq:fabrizio}) is shown to hold true by numerical
computations. Remark however that
\begin{eqnarray*}
\phi'(z,c) & = & c\cosh(cz) - \cosh(z) +1+c
\\
\phi''(z,c) & = & c^2\sinh(cz) - \sinh(z) \le0
\end{eqnarray*}
so $\phi'(0,c)=2c$, $\phi(\cdot, c)$ is concave,
$\lim_{z\to+\infty}\phi(z,c)=-\infty$ and $\phi(\cdot,c)$ has a
single point of maximum; at last $\phi'(z_c) = c\left(\cosh(cz_c)
-\cosh(z_c) \right)<0$. This concludes the proof of
(\ref{eq:eps1_c_alpha3}).

\smallskip

Remark that if $c=1/2$ then $\phi(z,c) = \sinh(z/2) - \sinh z +
\frac 32 z$. If $z$ is such that $\cosh z =3$, then $\sinh
z=2\sqrt 2$, $\sinh(z/2)=1$, $z=\log(3+\sqrt 8)$, and
(\ref{eq:fabrizio}) holds with strict inequality.
\end{remark}


\bigskip
\textbf{Acknowledgments.}  The authors kindly acknowledge support by
{\sc Progetto Gnampa 2005} \emph{Analisi Asintotica Per Si\-ste\-mi
  Iperbolici Non Lineari}.  The authors thank Graziano Guerra for
stimulating remarks and Umberto Massari for hints on $\BV$
functions.\small

\addcontentsline{toc}{section}{References} 

\end{document}